%% file: cohom1KarXiv4.tex
\documentclass[11pt]{article}
\usepackage{jeffplus_arxiv_c1K_b}
\numberwithin{equation}{section}
\usepackage{caption}
\usepackage[labelformat=simple]{subcaption}

\newextarrow{\xbigtoto}{{20}{20}{20}{20}}
{\bigRelbar\bigRelbar{\bigtwoarrowsleft\rightarrow\rightarrow}}
\newcommand{\longtoto}{\xbigtoto{}}

\newcommand{\niso}{\ncong}
\newcommand{\descend}{\vphantom{x_{X_{X_{X_{X_{X_{X}}}}}}}}
\renewcommand {\tE}{\wt{E}^* }
\let\hatorig\^

\renewcommand{\^}{\wedge}

\newcommand{\Og}{\mathsf{Orb}_G}

\newcommand{\man}{M}
\newcommand{\spac}{X}

\renewcommand{\X}{\Chi}
\renewcommand{\ring}{\Q}

\newcommand{\Km}{K^-}
\newcommand{\Kp}{K^+}
\newcommand{\Kpm}{K^\pm}
\newcommand{\apm}{\a^\pm}

\newcommand{\npm}{n_\pm}

\newcommand{\wpm}{w_\pm}

\renewcommand{\HH}{\H_{H}}

\newcommand{\Ke}{K_{\mr{eff}}}
\newcommand{\tKe}{\wt K_{\mr{eff}}}

\newcommand{\eff}{_{\mathrm{eff}}}
\newcommand{\He}{H_{\mr{eff}}}
\newcommand{\Hepm}{\He^\pm}
\newcommand{\tHe}{\wt H\eff}
\newcommand{\tHepm}{\wt H\eff^\pm}
\newcommand{\Kepm}{\Ke^\pm}
\newcommand{\tKepm}{\wt K\eff^\pm}

\newcommand{\MVS}{Mayer--Vietoris sequence\xspace}
\newcommand{\AHSS}{Atiyah--Hirzebruch spectral sequence\xspace}
\renewcommand{\ring}{k}
\newcommand{\casep}{join configuration\xspace}
\newcommand{\caseb}{sphere bundle configuration\xspace}

\begin{document}

\title{The K-theory of cohomogeneity-one actions}
\author{Jeffrey D.~Carlson}
\maketitle

\begin{abstract}
	\small{
	We compute the equivariant complex K-theory ring of
	a cohomogeneity-one action of a compact Lie group
	at the level of generators and relations
	and derive a characterization of K-theoretic 
	equivariant formality for these actions.
	Less explicit 
	expressions survive for a range of equivariant
	cohomology theories including 
	Bredon cohomology and Borel complex cobordism.
	The proof accordingly involves elements of
	equivariant homotopy theory,
	representation theory, and 
	Lie theory. 
	
	Aside from analysis of maps of representation rings
	and heavy use of the structure theory 
	of compact Lie groups,
	a more curious feature is the essential need 
	for a basic structural fact about the 
	Mayer--Vietoris sequence
	for any multiplicative cohomology theory 
	which seems to be otherwise unremarked in the literature,
	and a similarly unrecognized basic lemma 
	governing the equivariant cohomology 
	of the orbit space of a finite group action.
	}
\end{abstract}

\bs 

Compact Lie group actions $G \act M$ of \emph{cohomogeneity one},
those whose orbit space $M/G$ is a $1$-manifold,
have been a perennial object of study in differential geometry~%
\cite{mostert1957cohomogeneity,%
	neumann1968,%
	parker1986,%
	alekseevsky1993,%
	puettmann2009,%
	hoelscher2010classification,%
	frank2013cohomogeneity,%
	he2014obstruction,%
	galazgarciazarei2015,%
	angella2020cohomogeneity},
first because they are 
the 
most 
obvious 
class 
to study after 
homogeneous (= cohomogeneity-zero) actions,
but also 
because they furnish examples 
of Einstein metrics~\cite{berardbergery1982nouvelles} and 
manifolds with exceptional holonomy~\cite{bryantsalamon1989,cvetivc2002cohomogeneity,cvetivc2004new},
and especially because ``large'' isometry, 
for which low cohomogeneity gives a measure, 
has long played a central organizing role 
	(sometimes called the \emph{Grove program}~\cite{grove02ofandvia})
in 
finding Riemannian manifolds of nonnegative 
curvature~\cite{groveziller2000,%
	groveziller2002,%
	verdiani2004,%
	GVWZ2006,%
	grovewilkingziller2008,%
	ziller2009cohomogeneity,%
	dearricott2011,%
	verdianiziller2014}.
%
%
As nontrivial amounts of work have gone into understanding these actions geometrically,\footnote{\ 
	See the bibliography in the recent work of Galaz-Garc{\'i}a
	and Zarei~\cite{galazgarciazarei2015} 
	for some indication of the scope of this study.} 
their algebro-topological invariants are of some interest,
and phenomena arising in the computation of the 
rational Borel equivariant cohomology of these actions~\cite{CGHM2018}
hint at the generalization to a large class of cohomology theories
pursued in the present work.
The case of equivariant K-theory is particularly interesting, given 
its implications for the existence of vector bundles with prescribed properties;
for example, \Cref{thm:eqf} of the present work is used in a
work of Amann--Gonz\'alez-\'Alvaro--Zibrowius%
~\cite[Thm.~A(1)]{amann2019vector}
to construct metrics of non-negative curvature on vector bundles 
over a class of manifolds admitting cohomogeneity-one actions.

In considering cohomogeneity-one actions,
one almost always operates in the framework of
Mostert's classical structure theorem%
	\footnote{\ 
	with an important erratum caught 
	by Richardson and Samelson~\cite{mostert1957errata}
},
encapsulated in \Cref{fig:schematic}.
\begin{theorem}[Mostert%
~{\cite
	{mostert1957cohomogeneity}%
}]\label{thm:Mostert}
	Let $G$ be a compact Lie group acting smoothly on a compact
	smooth manifold $M$ in such a way that the quotient
	$M/G$ is a compact, connected $1$-manifold, possibly with boundary.\footnote{\ In the noncompact case, 
		where the quotient space is an open or half-open interval,
		$M$ deformation retracts onto a homogeneous fiber $G/H$
		of $M \lt M/G$,
		so this case is already understood from the point of view of this paper.}
\bitem
\item If $M/G$ is a closed interval, 
there are inclusions of closed subgroups 
$H \rightrightarrows \Kpm \rightrightarrows G$
such that $\Kpm/H$ are homeomorphic to spheres\footnote{\ 
	Without the smoothness hypothesis (omitted by Mostert), 
	$\Kpm/H$ can also be the Poincar\'{e} homology sphere,
	as noted by Galaz-Garc\'{i}a and Zarei only recently%
	~\cite
	{galazgarciazarei2015}.
	} and 
$M$ is
the double mapping cylinder 
of the span 
$G/H \rightrightarrows G/K^+$.
\item 
If $M/G$ is a circle,
there exist a closed subgroup $H$ of $G$
and an element $w$ of the normalizer $N_G(H)$
such that 
$M$ is diffeomorphic to the mapping torus of the right translation by $w$
on 
$G/H$.
\eitem
\end{theorem}	
\renewcommand\thefigure{\thesection.\arabic{figure}}  
\begin{figure}[H]
	\centering
	\begin{subfigure}[b]{5cm}	
		\centering
		\includegraphics[height=5.25cm]{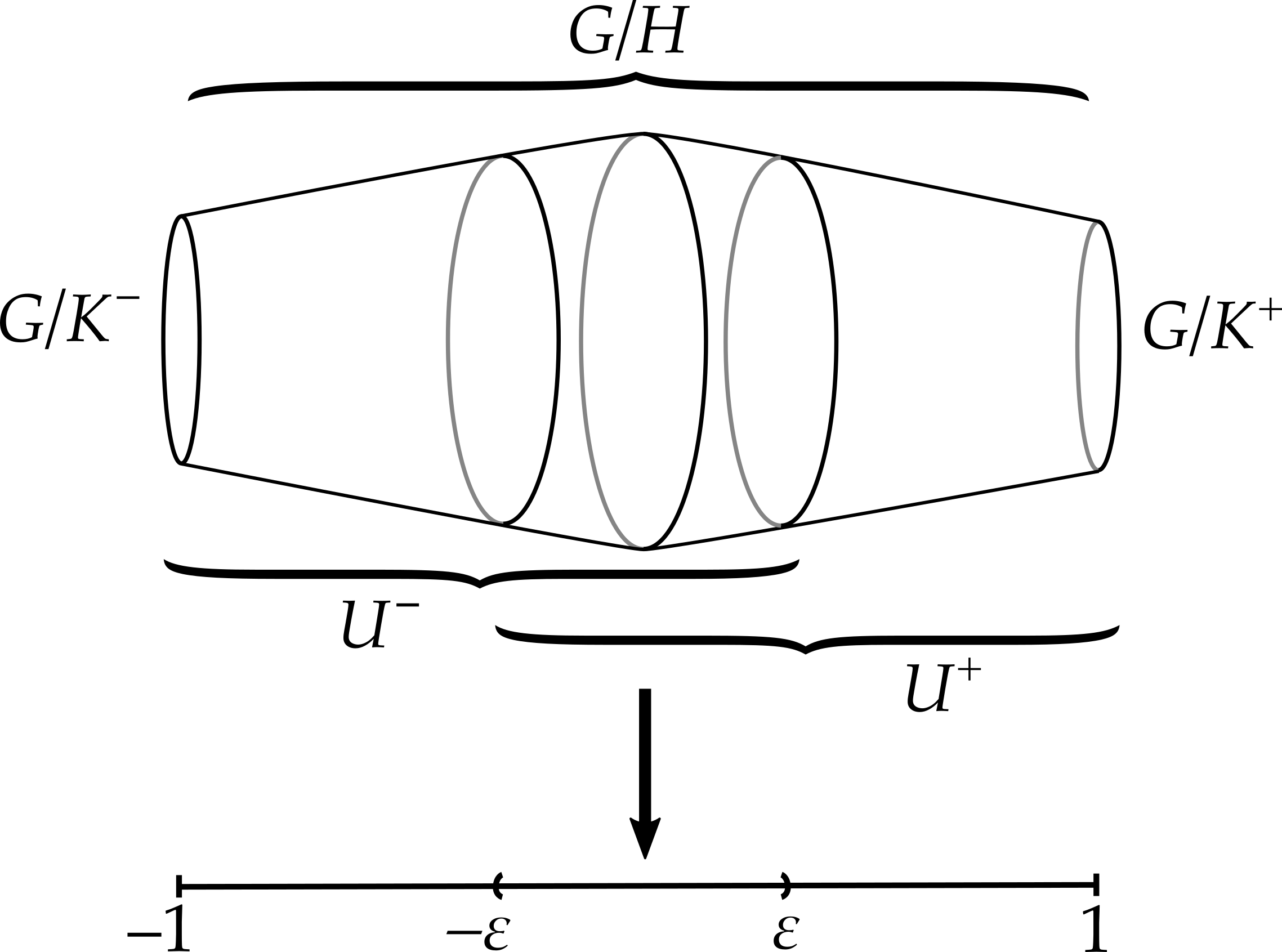}
		\caption{$M/G$ an interval}
		\label{fig:double}
	\end{subfigure}
	\hfill
	\begin{subfigure}[b]{5cm}
		\centering
		\includegraphics[height=5.25cm]{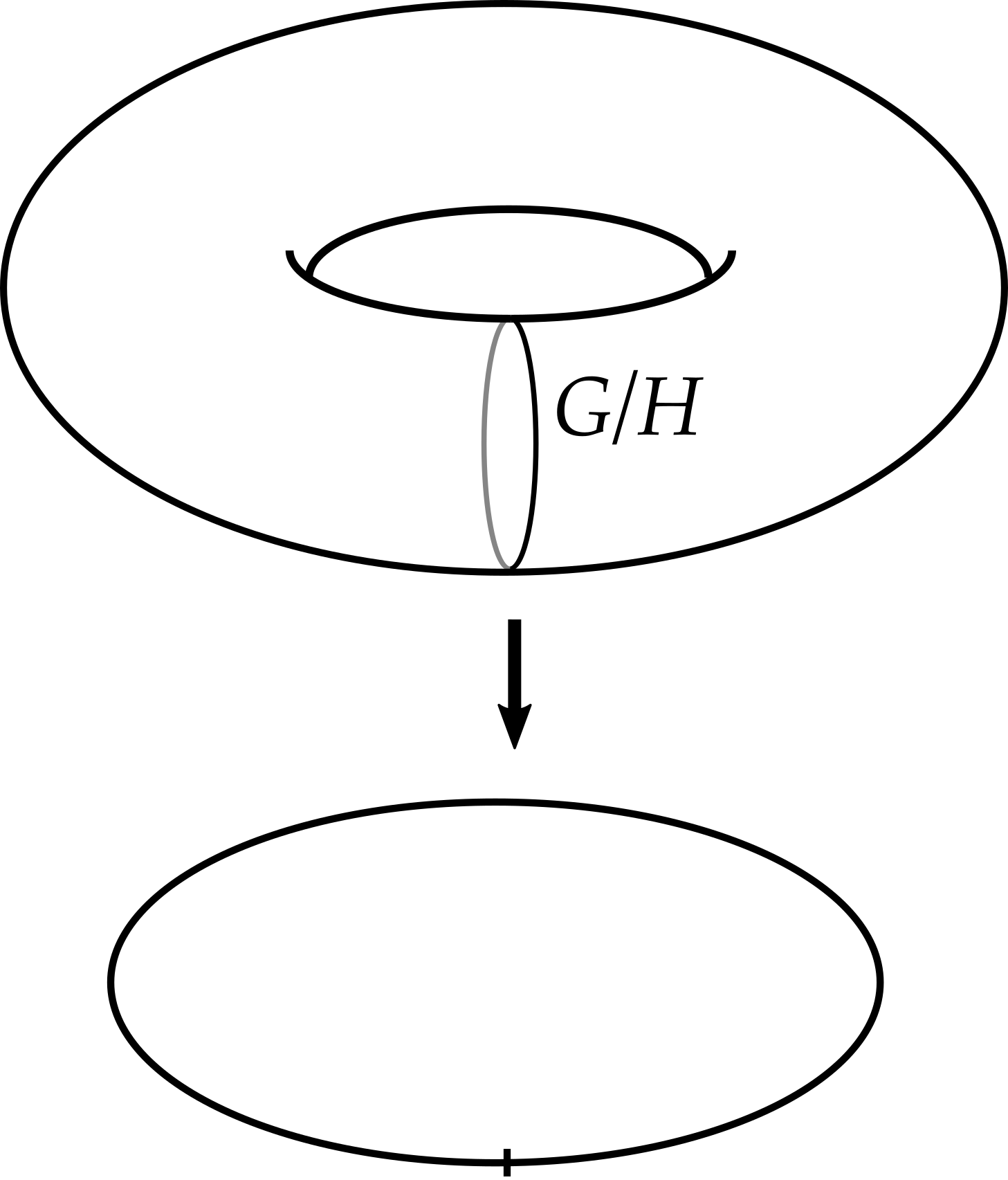}
		\caption{$M/G$ a circle}
		\label{fig:torus}
	\end{subfigure}
	\caption{Schematics for the orbit projection $M \lt M/G$ of a cohomogeneity-one action}
	\label{fig:schematic}
\end{figure}
%

In the case of the double mapping cylinder, if $M$ is smooth, 
then the isotropy quotients $\Kpm/H$ 
can actually be taken isometric in the Riemannian sense 
to round spheres
given by orbits in irreducible $\Kpm$-representations%
~\cite[Ex.~7.13]{Be},
suggesting equivariant complex K-theory $\defm{\KG}$,
whose coefficient ring is 
the ring $\defm{RG}$ of complex representations,
which is already motivated by its applications,
is also the most natural topological invariant of such an action. 
Indeed,
the \MVS of the cover $\{U^\pm\}$ in \Cref{fig:double} 
reduces to the exact sequence\quation{\label{eq:fourterm}
	0 
		\to 
	K^0_G (M) 
		\lt 
	RK^- \x RK^+ 
		\lt 
	RH 
		\os{\defm\d}\lt 
	K^1_G (M)
		\to 
	0\mathrlap,
}
where 
the middle map is the difference of the restrictions $R\Kpm \lt RH$
between complex representation rings,
showing the additive structure of $\KG(M)$ is wholly a question of representation theory.

Surprisingly, the multiplicative structure turns out to be as well.
The key fact 
is that the connecting map~$\d$ in \eqref{eq:fourterm} 
is actually a $K^0_G(M)$-module homomorphism.
The analogous fact in Borel cohomology
can be established by chasing cochains around a diagram,
but there are no cochains to follow in K-theory.
The result nevertheless turns out to be extremely general:

\begin{restatable*}{proposition}{MV}\label{thm:MV}
	Let $\E$ be a multiplicative ($\Z$-graded, $G$-equivariant) cohomology theory.
	Then the natural $\E(X)$-module 
	structure on the terms 
	of the Mayer--Vietoris sequence of a triad
	$(X;U,V)$ of $G$--CW complexes with $X = U \union V$
	is preserved by the connecting map
	in the sequence.
\end{restatable*}

\nd 
This basic result seems underappreciated;
working topologists surveyed by the author seem not to
know it,
nor does it seem to be discussed in the literature.
The enhanced connecting map makes life simpler in a variety
of situations, and a sample application to the cup product on a closed $3$-manifold is discussed in \Cref{ex:3mf}.
Most importantly for us, 
\Cref{thm:MV} immediately implies a general structure theorem 
for the equivariant cohomology ring of $G \act M$
in multiplicative cohomology theories with coefficients
concentrated in even degree, \Cref{thm:oddzero},
and one thus has a general expression for the
K-theory ring, \Cref{thm:mainK}.

To say more concretely what the ring $\KG(M)$ is,
one needs to explicitly identify 
the maps in the sequence \eqref{eq:fourterm}.
The structure theorem for $\HG(M;\Q)$
proceeds from analysis of an analogous sequence,
so one naturally changes the nouns in those statements
and hopes 
the same arguments will prove the stronger results.
While the results are indeed the expected ones,
the cohomological proof methods fail
utterly and the K-theoretic proof 
is incomparably more involved.

For example, the algebraic lemma
governing the map $\H(BK;\Q) \lt \H(BH;\Q)$ when $K/H$ 
is an odd-dimensional sphere 
is an easy result on commutative graded algebras,
but the analogous statement about surjections $RK \lt RH$ 
between ungraded polynomial rings
is a deep open problem in affine algebraic geometry,
the \emph{Abhyankar--Sathaye embedding conjecture}, 
and one is forced to an analysis in \Cref{sec:odd}
involving the structure theory of compact Lie groups
and the classification of homogeneous spheres.
The result when one of 
the spheres $\Kpm/H$ is odd-dimensional
then follows:

\begin{restatable*}{theorem}{Kodd}\label{thm:Kodd}
	Let $\man$ be the double mapping cylinder of the span 
	$G/H \longtoto G/\Kpm$ for inclusions
	$H \rightrightarrows \Kpm \rightrightarrows G$
	of closed, connected subgroups of a compact Lie group $G$
	such that $K^\pm/H$ are spheres 
	and the fundamental groups $\pi_1(\Kpm)$ are free abelian.
	
\medskip

		 (a) 
	Assume that both $\Kp/H$ and $\Km/H$ are odd-dimensional.
	Then we have an $RG$-algebra isomorphism of 
	$\KG(M) = K^0_G(M)$ with one of
	\[
	\frac
	{RH[t_-^{\pm 1},t_+^{\pm 1}]}
	{(t_- - 1)(t_+ - 1)}
	,\qquad\quad
	\frac
	{RH[t_-^{\pm 1},\ol\rho_+]}
	{(t_- -1)(\ol\rho_+)}
	,\qquad\quad
	\frac
	{RH[\ol\rho_-,t_+^{\pm 1}]}
	{(\ol\rho_-)(t_+  - 1)}
	,\qquad\quad
	\frac
	{RH[\ol\rho_-,\ol\rho_+]}
	{(\ol\rho_- \ol\rho_+)}\mathrlap,
	\]
	where we identify 
	$R\Kpm$ with the Laurent polynomial ring 
	$RH[t^{\pm 1}_\pm]$ when $\dim K^{\pm}/H = 1$
	and with the polynomial ring 
	$RH[\ol\rho_\pm]$ when $\dim K^{\pm}/H \geq 3$.
 	
\medskip
 	
 		(b) 
 	Assume 
 	$K^+/H$ is odd-dimensional
 	and $K^-/H$ is even-dimensional.
 	Then we have an $RG$-algebra isomorphism
 	of $\KG(M) = K^0_G(M)$ with 
 	\[
 	 R\Km \+ (t-1) RH[t^{\pm 1}] < RH[t^{\pm 1}] \iso R\Kp
 	\qquad\quad
 	\mbox{or}
 	\qquad\quad
 	 R\Km \+ \ol\rho RH[\ol\rho] < RH[\ol\rho] \iso R\Kp
 	 \mathrlap,
 	\]
	where we identify 
	$R\Kp$ with 
	$RH[t^{\pm 1}_\pm]$ if $\dim \Kp/H = 1$
	and with 
	$RH[\ol\rho_\pm]$ if $\dim \Kp/H \geq 3$.
 	 The product in either case is 
 	 determined by the restriction $R\Km \longmono RH$.

\medskip

	In all cases the $RG$-module structure is determined by restriction.
\end{restatable*}
 
Similar difficulties ensue when the spheres $\Kpm/H$
are both even-dimensional.
The determination of the product on $\HG(M;\Q)$ in this case 
reduces to pleasant arguments
involving {\SSS}s of fibrations between classifying spaces 
and the eigenspaces of the action of the so-called
\emph{Weyl group of a geodesic} of $M$ on $\H(BH;\C)$,
relying on the fact these eigenspaces 
are themselves graded vector spaces;
but the proof in K-theory 
involves a lengthy multi-layered induction
on the structure of compact Lie groups, 
whose base cases
require a number of lemmas in the Lie theory 
and representation theory
of simple Lie groups.
The result, however, comes out as clean as one could hope:
%

\begin{restatable*}{theorem}{Keven}\label{thm:Keven}
	Let $\man$ be the double mapping cylinder of the span 
	$G/H \longtoto G/\Kpm$ for inclusions
	$H \rightrightarrows \Kpm \rightrightarrows G$
	of compact Lie groups
	such that 
	the commutator subgroups
	of $\Kpm$ are products of simply-connected groups
	and $\SO(\odd)$ factors
	and $K^\pm/H$ are even-dimensional spheres. 
	Then there exist an element $z \in K^1_G(\man)$ 
	and an $RG$-algebra isomorphism
	\[
	\KG(\man) \,\iso\, 
	(RK^-|_H \inter RK^+|_H) 
	\,\ox\, 
	\ext[z],
	\]
	where the injections $R\Kpm \lt RH$
	and the $RG$-module structure are 
	given by restriction.
\end{restatable*}

\nd The base cases of the induction turn out to all be known special examples;
see \Cref{rmk:existingexamples}.

These structure results also allow one to characterize 
surjectivity of the map $\KG(M) \lt \K(M)$,
also known as K-theoretic \emph{equivariant formality},
using the Hodgkin--K\"unneth
and Atiyah--Hirzebruch--Leray--Serre spectral sequences
and some homological algebra:
\begin{restatable*}{theorem}{eqf}\label{thm:eqf}
	Consider a cohomogeneity-one action 
	of a compact, connected Lie group $G$
	with $\pi_1(G)$ torsion-free
	on a smooth closed manifold $M$
	such that the orbit space $M/G$ is an interval
	and the commutator subgroups
	of the exceptional isotropy groups $\Kpm$ 
	are the products of simply-connected groups
	and $\SO(\odd)$ factors.
	Then the action is K-theoretically equivariantly formal  
	if and only if $\rk G = \max\{\rk K^-,\rk K^+\}$.	
\end{restatable*}

So much for the case when $M/G$ is an interval.
When $M/G$ is a circle,
we can say nothing categorical 
before inverting the order $|\G|$ of the cyclic subgroup $\G$ 
generated by the class of $w \in N_G(H)$
in the component group $\pi_0 N_G(H)$ 
(see \Cref{ex:inversion}),
but once we do, the result follows formally 
from a much more fundamental fact about equivariant cohomology theories:

\begin{restatable*}{theorem}{covering}\label{thm:covering}
	Let $G$ be a compact Lie group
	and $\G$ a discrete finite group,
	and $X$ a finite $(G \x \G)$--CW complex
	whose isotropy subgroups are of the form 
	$H \x \D$ for $H \leq G$ and $\D \leq \G$.
	Moreover, let
	$\E$ be a $\Z$-graded $G$-equivariant cohomology theory
	valued in $\Z[\sfrac 1{\,|\G|}]$-modules.
	Then the quotient map $\pi\mnn\: X \lt X/\G$ induces an isomorphism
	\[
	\E(X/\G)  \isoto \E(X)^\G
	\] 
	onto the submodule 
	of $\G$-invariant elements.
\end{restatable*}
\nd The proof uses an equivariant \AHSS
and an observation about Bredon cohomology
to reduce to the classical result
for singular cohomology it generalizes,
and the result is again the sort of thing that one expects to find 
in the literature but does not.
In any event,
it has an immediate corollary, \Cref{thm:mappingtorus},
describing the equivariant cohomology of a mapping torus 
in broad generality,
which specializes to the result we wanted:
\begin{restatable*}{prop}{Kcircle}\label{thm:Kcircle}
Let $\man$
be the mapping torus
of the right translation by $w \in N_G(H)$
on a homogeneous space $G/H$
of a Lie group $G$ with finitely many components,
and write $\defm{w^*}$ 
for the maps induced on $\K(G/H)$
and $\KG(G/H) \iso RH$
by the right translation by $w$.
Let $\defm \ell$ be the least positive natural number
such that $w^\ell$ lies in the identity component of $N_G(H)$.
Then one has $\K (S^1)$- and 
$\big(RG \ox \K(S^1)\big)$-algebra isomorphisms
\eqn{
	\KG\big(\man;\Z[\sfrac 1{\,\ell}]\,\mn\big) 
		&\iso \xt{}
				{\K (S^1)}
				{(RH)^{\ang{r_w^*}}} \ox 
				\Z[\sfrac 1{\,\ell}]
						\mathrlap,\\
	\K\big(\man;\Z[\sfrac 1{\,\ell}]\,\mn\big) 
		&\iso \xt{}
				{\K (S^1)}
				{\K(G/H)^{\ang{r_w^*}}}\ox 
				\Z[\sfrac 1{\,\ell}]
						\mathrlap,
}
respectively, where 
$({-})^{\ang{w^*}}$ 
denotes the subring of $w^*$-invariant elements,
the $\K(S^1)$-module structure
is given in both cases by pullback from $M/G \homeo S^1$,
and the $RG$-algebra structure is induced 
by the inclusion $H \longinc G$.
\end{restatable*}

\medskip

The structure of the paper is as follows. 
The less involved case where $M/G$
is a circle, including \Cref{thm:Kcircle}, 
along with some necessary definitions, 
is discussed in \Cref{sec:mappingtori}.
In \Cref{sec:basic},
we assume the orbit space $M/G$ is an interval 
and discuss those
aspects of $\KG(M)$ which do not depend
on representation theory
on the dimensions of the homogeneous spheres $\Kpm/H$,
including the Mayer--Vietoris  \cref{thm:MV}
and a general structure \cref{thm:mainK}.
The refinements of this theorem
in the case $M/G$ is an interval,
depending on the parities of the dimensions of $\Kpm/H$,
rely on material on Weyl groups, Lie theory,
and maps of representation rings developed in \Cref{sec:mapRG}.
In \Cref{sec:odd}, 
we derive the consequences,
including \Cref{thm:Kodd}, when one of the spheres $\Kpm/H$
is odd-dimensional, and
in \Cref{sec:even}, we address the case when both of the spheres $\Kpm/H$
are even-dimensional and derive \Cref{thm:Keven}.
Finally, in \Cref{sec:eqf} we use these structural results 
to characterize K-theoretic equivariant formality 
for actions with orbit space an interval.

\medskip

\nd \emph{Acknowledgments.} 
The author would like to thank 
Omar Antol{\'i}n Camarena, 
Jason DeVito, 
Oliver Goertsches,
Chen He, 
Liviu Mare,
Clover May,
Marc Stephan, 
and
Marcus Zibrowius
for helpful conversations,
J\'an Min\'a\v{c} for thoughtful advice on presentation,
and the National Center for Theoretical Sciences in Taipei
for its hospitality during a phase of this work.

\input{tori}

\section{Mayer--Vietoris and double mapping cylinders}\label{sec:basic}

The circle case disposed of, 
we begin analyzing the double mapping cylinder \Cref{fig:double}
in Mostert's dichotomy \ref{thm:Mostert} from the introduction.

The double mapping cylinder $\man$ of $\pi^\pm\: G/H \longtoto G/\Kpm$
admits an obvious invariant open cover
by the respective inverse images $U^-$ and $U^+$
of the subintervals $[-1,\sfrac {1\mn} 2)$ and $(-\sfrac {1\mn} 2,1]$ of $\spac/G \homeo [-1,1]$,
and the intersection
$W = U^- \cap U^+$ equivariantly deformation retracts to $G/H$ and 
$U^\pm$ to $G/\Kpm$
in such a way that the inclusions $W \longinc U^\pm$
correspond to the projections $\pi^\pm$.
Since $\KG(G/\G) = K^0_G(G/\G) = R\G$ for closed subgroups $\G \leq G$
by restriction of an equivariant
bundle to the identity coset $1\G \in G/\G$ 
and $\KG$ is $\Z/2$-graded%
~\cite[Ex.~(ii), p.~132; Prop.~(3.5)]{segal1968equivariant},
the Mayer--Vietoris sequence in K-theory 
reduces to the exact sequence
\[
	0 
	\to 
	K^0_G (M) 
	\lt 
	RK^- \x RK^+ 
	\lt 
	RH 
	\os{\d}\lt 
	K^1_G (M)
	\to 
	0
\]
noted in the introduction.
As promised there,
this sequence is more informative than one might expect,
reflecting the fact that in great generality,
the properties of the \MVS are better than is commonly acknowledged.
Those who do not care about generality can safely 
substitute $\E = \KG$ everywhere in the following without loss.

%

\input{MVgeneral}

We can now finally return to K-theory.

\begin{restatable}{theorem}{mainK}\label{thm:mainK}
	Let $M$ be the double mapping cylinder
	of the projections 
	$\pi^\pm\: G/H \longtoto G/\Kpm$.
	The \MVS reduces to a short exact sequence
	\[
	0 \to K^0_G M \lt RK^-\mn \x RK^+ \lt RH \lt K^1_G M \to 0
	\]
	of $K^0_G M$-module homomorphisms,
	inducing the following
	graded ring and graded $RH$-module isomorphism,
	respectively:
	\[
	K^0_G (\spac) \,\iso\, \xu{RH}{R\Km\mn}{R\Kp}
	,\qquad
	K^1_G  (\spac) \,\iso\mspace{1mu} 
	\Big(\mn\quotientmed{RH\,}{\,R\Km|_H + R\Kp|_H}\Big)[1],
	\]
	where $(-)|_H$ denotes restriction of representations
	along the inclusions $H \longinc \Kpm$.
	The product of odd-degree elements is zero,
	and the product 
	$K^0_G (\spac) \x K^1_G (\spac) \lt K^1_G (\spac)$
	descends from the multiplication of $RH$:
	\[
	(\rho_-,\rho_+)\.\ol \s = \ol{\vphantom{X^{X^x}} \rho_-|_H \.\s }
	\]
	for $(\rho_-,\rho_+)$ in the fiber product $\xu{RH}{R\Km\mn}{R\Kp}$ 
	and 
	$\ol \s \in K^1_G(\spac)$
	the image of $\s \in RH$.
\end{restatable}


\bex
Let $G = \O(n)$ with $K = K^\pm = \O(3)$ and $H = \O(2)$ block-diagonal.
Recall that $R\O(3) \iso R\SO(3) \x R(\Z/2) =
\Z[\s,\e]/(\e^2-1)$,
where $\smash{\s\: \O(3) \inc \Aut \R^3 \to \Aut \C^3}$ 
complexifies the defining 
representation and $\e = \det\: \O(3) \lt \Aut \C$ is the determinant,
and $R\O(2) \iso \Z[\rho,\e]/(\e^2-1,\rho\e - \rho)$,
where $\rho\: \O(2) \lt \Aut \C^2$
complexifies the defining representation~\cite{minami1971representation}.
The restriction $RK \lt RH$
is given by $\s \mapsto \rho +1$ and $\e \mapsto \e$.
Now \Cref{thm:mainK} yields a short exact sequence
\[
0 \to K^0_G (\man) \lt \frac{\Z[\s,\e]}{(\e^2-1)} \x \frac{\Z[\s,\e]}{(\e^2-1)} \lt \frac{\Z[\rho,\e]}{(\e^2-1,\rho\e-\rho)} \to 0. 
\]
The kernel decomposes additively as the sum 
\[
\KG(\man) = K^0_G(\man) 
	\iso 
\Bigg\{(x,x): x \in \frac{\Z[\s,\e]}{(\e^2-1)}\Bigg\} 
\+ \big((\s-1)(\e-1),0\big) 
\+ \big(0,(\s-1)(\e-1)\big)
\]
This bears a familial similarity to the description in \Cref{thm:Keven}(b)
but cannot be put in those terms due to torsion.

The cohomological situation, by way of contrast, 
is much simpler: 
we have $\HK \iso \Q[p_1] \iso \HH$,
where $p_1$ the first Pontrjagin class of the tautological 
bundle over the infinite Grassmannian $\Gr(3,\Ri) = B\O(3)$,
so $\HG \man \iso \Q[p_1]$.
The equivariant Chern character 
taking a representation $V$ to the 
Chern character of the associated vector bundle $V_{\O(3)} \to B\O(3)$
sends $\s - 3$ to $p_1$ and annihilates $\e - 1$.
\eex

\bex
If $G = K^\pm = H$,  
the resulting double mapping cylinder is just the unreduced 
suspension $S(G/H)$ 
and
one has 
\eqn{
K_G^0 \big(S(G/H)\big) = RG, &\qquad\ 
K_G^1 \big(S(G/H)\big) \mspace{1.5mu} = \, \quotientmed{RH\,}{\,\im(\mspace{-1mu}RG \to RH)} \,[1].
}
\eex

\brmk\label{rmk:geom-gen}
The decomposition in \Cref{thm:mainK} 
admits a winning interpretation in terms of bundles. 
The isomorphism $R\Km\mn \x_{RH} R\Kp \isoto K_G^0(M)$ 
comes explicitly
from the decomposition of the double mapping cylinder as the union along $G/H$
of the mapping cylinders  
$\defm {M(}G/H \to G/\Kpm\defm)$
of the natural quotient maps $G/H \lt G/\Kpm$
for any pair $\s^\pm$ of $\Kpm$-representations agreeing on $H$,
one forms the union of the bundles 
$M(G \ox_H V_{\s^\pm} \to G \ox_{\Kpm} V_{\s^\pm}) \lt M(G/H \to G/\Kpm)$
along the restriction $G \ox_H V_{\s^\pm} \lt G/H$
to their common boundary.
Particularly, for a $\Kp$-representation $\s^+$ which is trivial on $H$,
one can extend the bundle $M(G \ox_H V_{\s^+} \to G \ox_{\Kp} V_{\s^+})$
by gluing on a trivial bundle over $M(G/H \to G/\Km)$; call this $\xi_{\s^+}$.
The formal difference of $\smash{\xi^+_{\s^+}}$ and the trivial bundle 
$\ul{\C}^{\dim V_{\s^+}}$ is a typical element of the summands 
$\ol\rho RH[\ol\rho]$ and  $(t-1) RH[t^{\pm 1}]$ figuring in \Cref{thm:Kodd}(a).

For \Cref{thm:Kodd}(b),
one similarly forms a virtual bundle $\xi_{\s^-}$ 
from a $\Km$-representation $\s^-$ trivial on $H$.
That the product $(\xi_{\s^-} - \ul\C^{\dim V_{\s^-}}) \ox 
(\xi_{\s^+} - \ul\C^{\dim V_{\s^+}})$ should be zero 
follows by noting the first factor is zero over $M(G/H \to G/\Km)$
and the second over $M(G/H \to G/\Kp)$. 

The map $RH \lt K^1_G(M)$ admits the following description.
Given an $H$-representation $\s$,
use Bott periodicity to send the class of the bundle $G \ox_H V_\s$
to an element of $K_G^0 \big(S^2 (G/H)\big)$,
and then pull back to an element of $K^0_G (SM)$ 
along the suspension of the map $M \longepi S(G/H)$
collapsing each of the end-caps $G/\Kpm$ to a point.

Hodgkin~\cite[Cor.~10.1]{hodgkin1975kunneth} notes 
the geometric significance of the class $\ol\b(\rho) \in K^1(K/H)$,
for $\rho$ a $K$-representation trivial on $H$, 
is as the class of the bundle on $S(K/H)$ 
obtained by gluing trivial bundles ${\ul{V\!}_\rho}$ 
over two copies of the cone $C(K/H)$ along 
their boundaries $K/H$ via the identification 
$(kH,v) \sim \big(kH,\rho(k)v\big)$.
\ermk

\section{Restrictions of representation rings}\label{sec:mapRG}

To say anything more meaningful about the map $RK^- \x RK^+ \lt RH$
figuring in \Cref{thm:mainK},
unsurprisingly,
we will have to do some representation theory.


\begin{definition}\label{def:balanced}\label{def:indecomposables}
If $\G$ is any group, 
we write $\defm{\G'}$ for its commutator
subgroup and $\defm{\G\ab}$ for its abelianization.
We then have a functorial short exact sequence $1 \to \G' \to \G \to \G\ab \to 1$.
The center of $\G$ 
is denoted by $\defm{Z(\G)}$ 
and the connected component of the identity element by $\defm{\G_0}$.
If two groups $\Pi$ and $A$
contain a subgroup $F$ central in both,
we write $\defm{\Pi \ox_F A}$
for the balanced product $(\Pi \x A) /  \big\{(f,f\-): f \in F\big\}$.
When a group $\G$ is isomorphic to such a balanced product
with $F$ finite, we refer to the isomorphism as a \defd{virtual product decomposition}.
It is well known that a compact, connected Lie group $\G$
admits a virtual product decomposition
$\smash{\dsp\G \iso \xt F {\G' }{Z(\G)_0}}$,
and $F$ is the intersection of $\G'$ and $Z(\G)_0$.

A representation ring $R\G$ is augmented over $\Z$
by the unique $\Z$-linear map taking an honest representation to its dimension.
Given a commutative ring $k$,
the category of augmentation-preserving maps of augmented $k$-algebras
is \emph{pointed} in the sense it admits $k$ as a zero object.
The kernel of the augmentation $A \lt k$ is denoted $\defm{IA}$,
or, if $A = R\G$ is a representation ring, $\defm{I\G}$.
The quotient $k$-module $IA/(IA)^2$,
the \defd{module of indecomposables}, is written $\defm{QA}$.
Specializing the general definition of exactness in a pointed category,
a sequence of augmented $k$-algebras $\smash{A \os f\to B \os g\to C}$
is said to be \defd{exact} at $B$ if 
$\ker g = f (IA) B$.	
A short exact sequence $k \to A \lt B \lt C \to k$
of augmented $k$-algebras is said to be \defd{split}
if there exists a section $C \longmono B$
inducing an isomorphism $A \ox_k C \isoto B$.

Given an inclusion $A \longinc B$ of rings, an element $b \in B$
is said to be \defd{transcendental} over $A$
if the $A$-algebra map $A[x] \lt B$ from
the polynomial ring in one indeterminate over $A$
sending $x$ to $b$ is injective.
\end{definition}

\subsection{The splitting lemma}

We need a refinement of the following splitting lemma due to Hodgkin.

\bthm[{\cite[Prop.~11.1]{hodgkin1975kunneth}}]\label{thm:Hodgkinsplit}
Given any compact, connected Lie group $K$ with free abelian fundamental group,
the sequence
\[
	\Z \to RK\ab \lt RK \lt RK' \to \Z
\]
induced by abelianization is split exact.
\ethm

This essentially allows us to factor out the representation ring
of the connected component of the center of a Lie group.
We actually want to factor out an arbitrary central torus.
In order for this to work we need $RK'$ to be a polynomial
ring, or equivalently,
that $K'$ be a direct product of 
simply-connected groups and odd special orthogonal groups~\cite{steinberg1975pittie}.\footnote{\ 
	The representation rings of the even special orthogonal
	groups and relation with those of the odd special 
	orthogonal groups are given in \eqref{eq:SO-surject}
	}

\bprop\label{thm:Hodgkinplus}
Let $K$ be a compact, connected Lie group 
such that $RK'$ is a polynomial ring
and let 
$\ul K$ be a connected subgroup containing $K'$ 
with free abelian fundamental group
and $A$ a virtual complement, 
meaning a central torus
with $\defm F = \ul K \inter A $
finite and such that $K \iso \ul K \ox_F A$.
Then the sequence
\[
\Z \to R(A/F) \lt RK \lt R\ul K \to \Z
\]
induced by the exact sequence $1 \to \ul K \to K \to A/F \to 1$
is split exact.
The splitting is not natural.

If $H$ is a closed, connected subgroup of $K$ 
also containing $A$
and $\ul H = \ul K \inter H$ contains $F$,
then the splittings can be chosen compatibly so that 
$RK \lt RH$
is identified with $R\ul K \ox R(A/F) \lt R\ul H \ox R(A/F)$
if 
\bitem
\item[(i)] the restriction $R\ul K \lt R \ul H$ 
is a split surjection or
\item[(ii)] the restriction $RK \lt RH$ is an injection.
\eitem
%
\eprop

\bpf
The proof of the first paragraph is the same as Hodgkin's,
but we reproduce it to verify it still works
if $\ul K$ is not semisimple.
From the assumption $RK'$ is a polynomial ring,
it follows from this lemma applied in the case
$\ul K = K'$ and $A = Z(K)_0$
that 
$R\ul K$ is the tensor product of a polynomial 
ring and a ring of Laurent polynomials.\footnote{\ 
	There is no circularity;
	one runs through this proof
	in the case $\ul K = K'$ is semisimple, 
	then once again in the general case.
	If this makes one uncomfortable, one can cite
	Hodgkin's \Cref{thm:Hodgkinsplit},
	noting that his proof never uses the assumption $K'$ 
	be simply-connected, 
	but just that $RK'$ is a polynomial ring.}
The restriction $\ul K \x A \lt K$ of the multiplication of $K$ 
is a surjective homomorphism with kernel the antidiagonal
 $\defm{\nu F} = \big\{(f,f\-): f \in F\big\}$ 
inducing the evident isomorphism $\ul K \ox_F A \isoto K$.
Pulling back, representations of $K$ can be identified 
with those representations of $\ul K \x A$ whose kernels contain $\nu F$.
The projections $\ul K \x A \epi A \epi A/F$
give us the first map $R(A/F) \longmono R(\ul K \x A)$
in the display.

For the second map, it will suffice to lift 
a list $(\defm{\rho_j})$ of representations of $\ul K$
forming a minimal set of polynomial and Laurent generators for $R\ul K$,
making sure the lifts of the Laurent generators are still units.
To lift an irreducible $\rho\: \ul K \lt \Aut \C^n$
to a representation of $\ul K \x A$ trivial on $\nu F$,
note that since $F$ is central, 
multiplication by each element of $\rho(F)$
is a $\ul K$-module endomorphism of $\C^n$,
and hence by Schur's lemma, a constant times $\id_{\C^n}$,
so $\rho|_F$ is a direct sum of $n$ copies of some
one-dimensional representation $\bar \s\: F \lt S^1$.
Since $\Hom(-,S^1)$ is exact and $F$ a subset of $A$,
taking $\rho = \rho_j$,
we see $\bar \s$ is the restriction of 
some $\defm{\s_{\rho_j}}\: A \lt S^1$.
For each $j$, consider the representation 
$\defm{\wt \rho_j} \ceq \rho_j \ox (\s_{\rho_j})^{\+n}$ of $\ul K \x A$
in $\C^n$
taking $(k,a) \lmt \s_{\rho_j}(a)\id_{\C^n}\.\rho_j(k)$.
This $\wt \rho_j$ vanishes on $\nu F$ by construction and restricts to $\rho_j$ on $\ul K$.
In case $\rho_j\: \ul K \lt S^1$ was one of the Laurent generators,
then $n=1$, so $\wt\rho_j$ is still a one-dimensional representation
and hence invertible.

It remains to show the map is an isomorphism. We have maps
\[
	R\ul K \ox R(A/F) \os{\defm{\phi}}\lt RK \longmono R\ul K \ox RA,
\]
where $\phi$ is defined in the expected manner 
from the maps we have just constructed
and the second map comes from the covering $\ul K \x A \lt K$
and the natural identification $R(\ul K \x A) \iso R\ul K \ox RA$.
Since $A \longto A/F$ is surjective,
$\Hom(A/F,S^1) \lt \Hom(A,S^1)$
and hence $R(A/F) \lt RA$ are injective.
Hence the composition is injective
on elements of the form $p(\vec \rho) \ox \t$,
where $p(\vec \rho)$ is a Laurent monomial in the generators $\rho_j$
and $\t$ is an element of $\Hom(A/F,S^1)$.
As such elements form a $\Z$-basis for 
$R\ul K \ox R(A/F)$,
we find $\phi$ is injective.
To see it is surjective, 
let any element $p(\rho_j) \ox \t \in R(\ul K \x A)$
vanishing on $\nu F$ be given; such elements
form a $\Z$-basis for the image of $RK \longmono R\ul K \ox RA$.
The element can be rewritten 
$p(\rho_j) \ox \t = p({\wt\rho}_j) \. (1 \ox \t')$
for some other $\t' \in \Hom(A,S^1)$.
Moreover, $1 \ox \t'\: (k,a) \lmt \t'(a)$ is trivial on $\nu F$ 
since $p(\rho) \ox \t$ and $ p(\wt{\rho}_j)$ are,
so $\t'$ is trivial on $F$ and hence descends to an element of $R(A/F)$.
Thus $p(\rho_j) \ox \t = \phi\big(p(\rho_j) \ox \t'\big)$.

\smallskip

Now we bring in $H$. 

\smallskip

\nd (i)
There is a natural map from $1 \to \ul H \to H \to A/F \to 1$
to the exact sequence for $K$,
inducing a map of short exact sequences
of representation rings.
A choice of splitting $R\ul H \to R\ul K$
and the splitting $R\ul K \to RK$ of the first part of the proposition
uniquely induces a 
compatible splitting
$R\ul H \to R\ul K \to RK \to RH$.

\smallskip

\nd (ii)
We have already seen that $RK \lt R\ul K$
and $RH \lt R\ul H$ are surjections
with kernels generated by the image of $I(A/F)$.
It follows $R\ul K \lt R\ul H$ is a monomorphism
which can be identified with the reduction of $RK \lt RH$
modulo $I(A/F)$.
Let a splitting $\phi$ of $RH \longepi R\ul H$
be given, and consider the image $R$ of
the composition $R\ul K \mono R\ul H \mono RH$. 
By definition, the subring generated by 
$R$ and the image of the natural map $R(A/F) \lt RH$
is abstractly isomorphic to $R\ul K \ox R(A/F)$
and surjects onto the image of $R\ul K \longmono R\ul H$,
so this subring is the image of $RK \lt RH$.
\epf

\subsection{Lemmas for odd spheres}

The results we need for the case 
the homogeneous sphere $K/H$ is odd-dimensional all follow 
from the splitting \cref{thm:Hodgkinplus} 
once we show $RK \lt RH$ is split surjective.

\bprop\label{thm:Koddcircle}
Let $H \leq K$ be connected, compact Lie groups 
such that $K/H \homeo S^1$ and $RK'$ is a polynomial ring.
Then $RK \lt RH$ is a surjection
and can be written 
\[
	RH[t^{\pm 1}]
		\xrightarrow{t \mapsto 1} 
	RH,
\]
where 
$t\: K\ab \epi K\ab/H\ab \simto \U(1)$ 
pulls back one of the generators of $R(K\ab/H\ab)$
and is transcendental over $RH$.
\eprop
\bpf
Consider the diagram
\quation{\label{eq:abelianization}
	\begin{aligned}
	\xymatrix@C=1em{
		1\ar[r]&H' \ar[r]\ar[d]	& H \ar[r]\ar[d]	& H\ab \ar[d]\ar[r]& 0\\
		1\ar[r]&K' \ar[r]	& K \ar[r]\ar[d]	& K\ab \ar[d]\ar[r]& 0\\
			   &			& K/H \ar[r]		& K/K'H, 
		}
	\end{aligned}
}
whose first two rows are exact sequences and 
whose first two rows and second column are fibrations.
Since $\pi_2$ of a Lie group is zero,
and $\pi_1 (H')$ and $\pi_1 (K')$ are finite,
we see $\pi_1(K/H) \ox \Q \to \pi_1(K/K'H) \ox \Q$ 
is an isomorphism,
so the torus $K\ab/H\ab = K/K'H$ is a circle.
Particularly, it is one-dimensional,
so counting other dimensions, 
we have
$\dim K' = \dim H'$, 
meaning $K'/H'$ is a connected $0$-manifold
and hence $K' = H'$.

The exact sequences of representation rings 
resulting from the first two rows of (\ref{eq:abelianization}) split
by \Cref{thm:Hodgkinplus}.
These splittings are not natural, but since $RK' \lt RH'$ is an isomorphism,
we can choose the liftings compatibly so that the following diagram commutes:
\quation{\label{eq:repsplit}
	\begin{aligned}
	\xymatrix@C=.0625em@R=.5625em{ 
	&&RK\ab \ox RK'\ar[dd]^\vertsim &&\qquad \	\\
	\\
	 \Z \ar[r] & RK\ab\!\! \ar[r]\ar[ruu]\ar[ddd] & RK \ar[r]\ar[ddd] & RK' \ar[r] \ar[luu]\ar[ddd]^\vertsim& \Z\phantom{.}\\
	 &&&\\
	 &&&\\
	 \Z \ar[r] & RH\ab\!\! \ar[rdd]\ar[r] & RH \ar[r] & RH'\ar[dld] \ar[r] & \Z.\\
	 \\
\qquad \ 	&& RH\ab \ox RH'\ar[uu]_\vertsim&&
	}
	\end{aligned}
}
Since $RK\ab \lt RH\ab$ is induced by the inclusion 
$H\ab \longinc K\ab$ of a codimension-$1$ subtorus 
and monomorphisms between tori admit retractions, 
we have $RK\ab \iso RH\ab \ox R(K\ab/H\ab) \iso RH\ab[t^{\pm 1}]$
and the result follows.
\epf

\bprop\label{thm:Koddsphere}
Let $H \leq K$ be connected, compact Lie groups 
such that $K/H$ is a sphere of odd dimension $3$ or more
and $RK'$ is a polynomial ring.
Then $RK \lt RH$ is a surjection
and if a subgroup of $K$
surjecting on the image of $K \lt \Homeo K/H$ 
is simply-connected, then $RK \lt RH$ can be written as 
\[
	RH[\bar\rho] 
		\xrightarrow{\bar \rho \mapsto 0} 
	RH,
\]
where $\bar \rho$ is transcendental over $RH$ and 
equals $\rho - \dim \rho$ 
for a $K$-representation $\rho$, trivial on $H$,
such that the induced continuous map $K/H \lt \U$ 
represents the fundamental class of $K/H$.
\eprop
\bpf
In (\ref{eq:abelianization})
the bottom map now is a fibering of an odd sphere
over a torus, which is only possible
if the torus in question is zero-dimensional.
Thus $H\ab \lt K\ab$ is a homeomorphism, 
so $H' = \ker(H \to H\ab)$ and 
$H \inter K' = \ker(H \to K \to K\ab)$ are equal.
Since $K/K'H$ is trivial
and the fiber of the trivial map $K/H \lt K/K'H$
is $K'H/H \iso K'/(K' \inter H) = K'/H'$,
it follows $K'/H' \lt K/H$ is a homeomorphism.
By the following \Cref{thm:spheresurj}, 
one has $RK' \lt RH'$ a surjection
of the form 
$RH'[\bar\rho] \simto RK' \to RH'$
if the group $\tKe'$ of that lemma can be taken simply-connected,
so \Cref{thm:Hodgkinplus}(i)
applies with $A$ the maximal central torus of $H$
and $\ul K = K'$ and $\ul H = H'$.


To show the generator has claimed property, 
recall that the Hodgkin map $\b\: R\G \lt \K(\G)$ is functorial,
factors through the 
{module of indecomposables} ${QR\G}$,
and induces isomorphisms
$\ext_\Z [QR\G] \isoto \K(\G)$
if $\pi_1 (\G)$ 
is torsion-free,
as we now assume $\pi_1(K)$ (and hence $\pi_1(H)$) is.
Thus $i^*\: \K(K) \lt \K(H)$
is a surjection.
A result of Minami~\cite[Prop.~4.1]{minami1975symmetric}
then says $\K(K/H)$ is the exterior algebra on
the homotopy class $\ol\beta(\rho)$ of 
the composition $K/H \to \U(V_\rho) \inc \U$
for an element $\rho \in RK$ 
whose class in $QRK$ generates $\ker Q(RK \to RH)$.
\epf

We have separated out the harder part of the preceding proof 
into that of the following result.

\bprop\label{thm:spheresurj}
Let $H \leq K$ be connected, compact Lie groups 
such that $K/H$ is a sphere of odd dimension $3$ or more
and $RK'$ is a polynomial ring.
Then the map $RK' \lt RH'$
is an augmentation-preserving surjection
which can be written as $RH'[\bar\rho] \longepi RH'$
for a judicious choice of section $RH' \longmono RK'$ and 
algebraically independent generator $\bar\rho$.
\eprop
\bpf
Recall $K'$ is a direct product of 
simply-connected simple groups and odd special orthogonal groups~\cite{steinberg1975pittie}.
The action of $K'$ on $K'/H' \homeo S^{2n-1}$
induces a homomorphism $K' \lt \Homeo S^{2n-1}$
whose image we dub $\defm{\Ke'}$.
The induced surjection $\f k' \lt \f k'\eff$
of semisimple Lie algebras splits.
The exact sequence of the fibration $K' \to K \to K\ab$ implies 
the finite group $\pi_1 K'$ vanishes,
so there is some subgroup
$\defm{\wt K'\eff}$ of $K'$
such that the composition 
$\wt K'\eff \inc K' \epi \Ke'$ is a finite covering.
The structure theorem for connected, compact Lie groups thus implies
$\wt K'\eff$ is a direct factor of $K'$, say
$\smash{K \iso \wt K'\eff \x L}$.
Note that $L$ lies in the kernel of $K' \longepi K'\eff$,
and so particularly is contained within $H'$,
so that
$H' \iso \wt H'\eff \x L$,
where $\defm{\wt H'\eff} \ceq H' \cap \wt K'\eff$.
Thus we may write $RK' \lt RH'$
as $\id_{RL} \ox \,(R\tKe' \to R\tHe')$
and we need only analyze the last factor.
Augmentation-preservation is just the fact 
restriction of representations preserves dimension,
so it remains only to see $R\wt K'\eff \lt R\wt H'\eff$
is a surjection of the claimed form.
This comes down to a short case analysis,
%
%
as the entire list of expressions for an odd-dimensional sphere
as the orbit of an effective action of a compact, connected Lie group
is the following~\cite[Ex.~7.13]{Be}\cite[Table~C, p.~104]{grovewilkingziller2008},
where the balanced product notation $\ox_{\Z/2}$
is as explained in \Cref{def:balanced}:
\quation{\label{eq:spherelist}
	\begin{aligned}
		S^{4n-1} &= \frac{\Sp(n)}{ \Sp(n-1)} = 
		\frac{\xt{\Z/2}{\U(1)}{\Sp(n)}\descend}
		{\xt{\Z/2}{\Delta \U(1)}{\Sp(n-1)}\descend} =\frac{\xt{\Z/2}{\Sp(1)}{\Sp(n)}\descend}
		{\xt{\Z/2}{\Delta \Sp(1)}{\Sp(n-1)}\descend} 
		,\\
		S^{2n-1} &= \U(n) /\U(n-1) = \SU(n) /\SU(n-1) = \SO(2n) /\SO(2n-1),\\
		S^{15} &= \Spin(9) / \Spin(7),\\
		S^7 &= \Spin(7) / G_2.
	\end{aligned}
}
Our task is made easier by the $\l$-ring structure on $R(-)$
induced by exterior powers:
because the rings in question are largely generated
by exterior powers of the standard representation $\defm \s$,
much of the work is done when we find $\s$ in the image.
\bitem
\item
For $R\Sp(n) \lt R\Sp(n-1)$ 
we have $\s \lmt \s + 2$
and for
$R\SU(n) \lt R\SU(n-1)$ 
we have $\s \lmt \s + 1$.
Now $\s$ generates $R\Sp(n)$ and $R\SU(n)$ as $\l$-rings,
so we already see the map is surjective.

In fact, the images of $\s,\ldots,\l^{n-1}\s$ 
generate the codomain in either case,
since $\l^j(\s+2) = \l^j \s + 2 \l^{j-1}\s + 1$ for $j \geq 2$ 
and $\l^j(\s+1) = \l^j \s + \l^{j-1} \s$ for $j \geq 1$.\footnote{\ 
	In general $\l^n(x+y) = \sum_{i+j = n} \l^i x \. \l^j y$,
	and for $m \in \N$ one has $\l^j m = {m \choose j}$.
}
It follows the image of $\l^n \s$ is also the image of 
some polynomial $p$ in the lower $\l^j \s$,
so we may rewrite the domain as $\Z[\s,\ldots,\l^{n-1}\s][\l^n \s - p]$
to obtain an expression of the claimed form.
\item
Writing $R\Spin(2n) \lt R\Spin(2n-1)$
as 
$\Z[\s,\ldots,\l^{n-2}\s,\D_-,\D_+] \lt \Z[\s,\ldots, \l^{n-2}\s,\D]$,
where $\s$ is the composition of the double cover with the defining representation of the special orthogonal group,
$\defm{\D_\pm}$ are the half-spin representations,
and $\defm\D$ is the spin representation,
we have $\s \lmt \s + 1$ and $\D_\pm \lmt \D$%
~\cite[Prop.~VI.6.1]{brockertomdieck}.

By the same argument as before, the map is a bijection
when restricted to $\Z[\s,\ldots,\l^{n-1}\s,\D_-]$,
and we may replace the last generator by $\D_+ - \D_-$
to obtain the desired expression.
\item 
The restriction $R\SO(2n) \lt R\SO(2n-1)$
is surjective because representations of $\SO(2n-1)$
descend from representations of $\Spin(2n-1)$
such that $-1 \in \Spin(2n-1)$ acts trivially,
and we have just seen the map $R\Spin(2n) \lt R\Spin(2n-1)$
is surjective.\footnote{\ 
	We will not use this case further, as $\SO(2n)$ is not simply-connected,
	but it is worth laying out clearly.
}

To get more specific expressions, we~\cite[Prop.~VI.6.6]{brockertomdieck}
may write the map as
\quation{\label{eq:SO-surject}
	\Z[\s,\ldots,\l^{n-1}\s,\l^n_+,\l^n_-]/(Q) 
		\lt 
	\Z[\s,\ldots,\l^{n-1}\s]\mathrlap,
}
where $\defm{\l^n_\pm}$ are the $\pm 1$-eigenspaces of the Hodge star on 
$\l^n \s$
and 
\[
\defm Q = (\overbrace{\l^n_+ + \l^{n-2}\s + \cdots}^{\defm x})
	(\overbrace{\l^n_- + \l^{n-2}\s + \cdots}^{\defm y})
		- (\overbrace{\l^{n-1}\s + \l^{n-3}\s + \cdots}^{\defm z})^2.
\]
We have a decomposition $\l^n \s = \l^n_+ + \l^n_-$ into
irreducibles, 
and $\l^n \s \lmt \l^n \s + \l^{n-1} \s = (\l^{n-1} \s)\dual + 
\l^{n-1} \s = 2\l^{n-1}\s$
in $R\SO(2n-1)$ since the fundamental representations of $\SO(2n-1)$
are self-dual, 
so it follows that both of $\l^n_\pm$ are sent to $\l^{n-1}\s$.
If we rewrite $R\SO(2n)$ as $\Z[\s,\ldots,\l^{n-2}\s][x,y,z]/(xy-z^2)$,
we see each of $x,y,z$ map to $w = \sum_{j \leq n-1} \l^j\s$
(so particularly $Q \lmt 0$),
and the map can finally be viewed as
\quation{\label{eq:SOmap}
	\quotientmed{\Z[\s,\ldots,\l^{n-2}\s][x,y,z]\,}{(xy-z^2)}
		\,\lt\,  
	\Z\big[\s+1,\ldots,\l^{n-2}(\s+1)\big][w].\phantom{,z]}
}

\item 
One~\cite{LiEonline}
can write $R\Spin(7) \lt RG_2$ as 
\eqn{
\Z[\s,\l^2\s,\d] 	&\lt \Z[\s,\Ad],\\ 
				\s 	&\lmt \s,\\ 
				\d 	&\lmt 1 + \s,\\ 
	\Ad = \l^2\s 	&\lmt \l^2\s = \s + \Ad.
}
Particularly, one can obtain the desired expression
by
exchanging the generator $\l^2\s$ for $\l^2\s-\s$ and $\d$ for $\d-\s-1$.
\item
One~\cite{verdianiziller2009,LiEonline} can write
$R\Spin(9) \lt R\Spin(7)$
as 
\eqn{\Z[\s,\l^2\s,\l^3\s,\D] &\lt \Z[\s,\l^2\s,\d],\phantom{\l^3 \s,}\\
	\s &\lmt \d + 1,\\
	\D &\lmt \d + \s + 1.
}
Then we have $\l^2(\s-1) \lmt \l^2 \d = \s + \l^2 \s$
and $\l^3(\s-1) \lmt \s\d-\d$.
Thus we can take instead as generators 
\eqn{
\s -1 &\lmt  \d,\\
\D - \s &\lmt \s,\\
\l^2(\s-1) - (\D-\s) &\lmt \l^2\s,\\
\l^3(\s-1) - (\D-\s-1)(\s-1) &\lmt 0.\qedhere\phantom{\l^3(\s-1) - (\d-\s-1)(\s)}
}

\eitem
\epf

\brmk 
The two ``exceptional'' homogeneous spheres can be
understood as follows.
Recall that the compact exceptional group $G_2$ can be seen as the group
of $\R$-algebra automorphisms of the octonions $\Oct$. 
The map $G_2 \longmono \Spin(7)$ lifts the inclusion 
$G_2 \longinc \SO(7)$ arising from restriction of the defining action
to the subspace of pure imaginaries. 
For the map $\Spin(7) \longmono \Spin(9)$ yielding $S^{15}$,
since $\pi_1\Spin(7) = 1$, one lifts the spin representation 
$\d\:\Spin(7) \longmono \SO(8)$
to $\Spin(7) \longmono \Spin(8)$,
then follows with the map $\Spin(8) \longinc \Spin(9)$
double-covering the block-diagonal inclusion $\SO(8) \+ [1] \longinc \SO(9)$.

The author learned these explanations from Jason DeVito.
\ermk

\brmk\label{thm:polysurj}
The proof of \Cref{thm:spheresurj}
was originally routed through the following statement:

\medskip

\emph{For any 
	surjection 
	$\varphi\: A \longepi B$
	of polynomial rings respectively
	in $m \geq n$ indeterminates 
	over a commutative base ring $\ring$, 
	one can choose an algebraically independent 
	set $x_1,\ldots,x_n,y_{n+1},\ldots,y_m$
	of polynomial generators for $A$ over $\ring$
	such that $\varphi$ sends $y_j \lmt 0$ 
	and restricts to an isomorphism $\ring[x_1,\ldots,x_n]\isoto B$.
}

\medskip

This innocuous-sounding claim 
is true for graded maps of graded rings over $\ring = \Q$
and open for ungraded maps over $\ring = \C$.
In algebro-geometric language, 
the special case $m = n+1$ we use in this paper
is the \emph{Abhyankar--Sathaye embedding conjecture}%
~\cite{abhyankarmoh1975,%
	sathaye1976linear,%
	russellsathaye,%
	popov2015around,%
	MO:Wendt},
which states that any embedding 
$\mathbb{A}^{\mn n}_\C \longmono \mathbb{A}_\C^{\mn n+1}$
is taken to the standard embedding 
by some automorphism of $\mathbb{A}_\C^{\mn n+1}$.
This is known at present for $n = 1$ and several other special cases, 
and is closely related to the determination of the algebraic automorphism group
$\Aut \mathbb{A}_\C^{\mn m}$, which is still incomplete for $m \geq 3$.
\ermk

\subsection{Lemmas for even spheres}\label{sec:even-lemma}

In case the homogeneous sphere $K/H$ is even-dimensional, 
the restriction $RK \lt RH$
makes the $RH$ a free module of rank two over $RK$.

\bprop\label{thm:Kevensphere}
Let $H \leq K$ be connected, compact Lie groups of equal rank
such that $K/H$ is an
even-dimensional sphere
and the semisimple component $K'$ 
is the direct product of a simply-connected group
and $\SO(\odd)$ factors.
Then $RH$ is a free $RK$-module of rank two.
\eprop
\bpf
Steinberg~\cite{steinberg1975pittie},
strengthening an earlier result of Pittie,
shows
that with our hypotheses, $RH$ is free of rank $|W_K|/|W_H|$ over $RK$ (he also provides a basis). 
To see the rank is two, 
note that by completion~\cite[Thm.~5.3]{carlsonfok2018},
this is also the rank of
$\H(BH;\Q)$ over $\H(BK;\Q)$,
which is $2$ by the collapse of the 
\SSS of $K/H \to BH \to BK$ with rational coefficients.
\epf

To apply this we will often use the splitting in
\Cref{thm:Hodgkinplus}(ii),
and for this we need to check that the condition on the finite
subgroup $F$ is satisfied.

\blem\label{thm:lowerF}
Suppose a compact, connected Lie group $K$ 
can be written as balanced product $\ul K \ox_F A$
of two subgroups $A$ and $\ul K$,
where $A$ is a central torus in $K$ and $F$ is finite,
and that $H$ is a closed subgroup of $K$ such that
$K/H$ is a sphere $S^{2n}$ of positive even dimension.
Then, writing $\ul H = H \inter \ul K$,
we have $H \iso \smash{\xt F{\ul H}A}$,
where the balanced product notation is recalled in  \Cref{def:balanced}.
\elem
\bpf
Since $\pi_1 (S^{2n}) = 0$,
it follows $H$ must contain $A$,
and it follows from the decomposition of $K$ 
that $\ul H$ and $A$ together generate $H$.
The preimage of $H$ under the projection 
$\ul K \x A \lt K$ is $F \ul H \x A$,
so it follows $\ul K / F\ul H \homeo S^{2n}$.
Since $\ul K/\ul H \lt \ul K/F\ul H$ is a finite covering, 
we see $F \ul H = \ul H$,
so $\ul H$ contains $F$.
Thus one can write $H \iso \ul H \ox_F A$ as claimed.
\epf

We will use this reduction in conjunction with a refinement
due to Adem and G{\'o}mez
of the Steinberg basis theorem.

\bthm[{Adem--G{\'o}mez \cite[Thm.~3.5]{ademgomez2012}}]\label{thm:AGbasis}
Let $G$ be a compact, connected Lie group with free abelian fundamental group and
fix a choice $\Phi^+$ of positive roots of $G$
with respect to some maximal torus.
Let $\ms W = (W_j)$ be a family of subgroups of $W = WG$, 
including the trivial group $1$
and $W$ itself, each generated by reflections in some subsystem $\Phi_j^+$
of $\Phi^+$,
and write
$\smash{\wt{W/W_j}} \ceq \{w \in W : w\Phi^+_j \sub \Phi^+\}$
for each $j$; this is a set of coset representatives for $W/W_j$.
Suppose any pair $W_j$ and $W_k$ in $\ms W$
lie in a common supergroup $W_\ell \in \ms W$ such that 
$\wt{W/W_\ell} = \wt{W/W_j} \inter \wt{W/W_k}$. 
Then $RT$ is the free $(RT)^W$-module $(RT)^W \. B$
on a basis $B \sub RT$
in such a way that each $j$ is associated to a subbasis $B_j \sub B$
such that $(RT)^{W_j} = (RT)^W \. B_j$ under the identification
and further, for any containment $W_k \geq W_j$ in $\ms W$,
there is a corresponding containment $B_k \sub B_j$
such that the induced inclusion
$(RT)^W \. B_k \longinc (RT)^W \.B_j$
is identified with the inclusion of invariants $(RT)^{W_k} \longinc (RT)^{W_j}$.
\ethm

We will apply this lemma in a number of cases,
invoking some elementary facts about extensions of root systems.

\blem\label{thm:latticeinclusionG2}
A lattice of Killing--Cartan type $A_2$
extends to a $G_2$ lattice in a unique way.
\elem
\bpf
If view the $A_2$ lattice
as the vectors $(a_1,a_2,a_3) \in \Z^3$
with $a_1 + a_2 + a_3 = 0$,
a new simple root $\a$ in an extending $G_2$ lattice
must have length $\sqrt 6$ 
and inner products with $A_2$ lattice elements divisible by $3$.
We would not rob the reader of the simple joy of verifying
only $\pm(2,-1,-1)$, $\pm(-1,2,-1)$, and $\pm(-1,-1,2)$ do the job.
\epf

\blem\label{thm:latticeinclusionBn}
A lattice of Killing--Cartan type $D_n$ 
extends to a $B_n$ lattice
in\\
\[
\case{ \mbox{a unique way } &\mbox{if }n \neq 4,\\ 
		\mbox{precisely two ways} & \mbox{if } n = 4.
		\vphantom{X_{X_{X_{X_x}}}}}
\]

	
\elem
\bpf
The standard $D_n$ lattice in $\R^n$ 
is spanned by roots $e_j \pm e_k$,
and so is given by those
integer linear combinations 
$\sum a_j e_j$ of the standard basis vectors $e_j \in \R^n$ 
for which $\sum a_j$ is even.
A new root $\a$ in an extending $B_n$ lattice 
must have length $1$ and inner product 
with all such vectors an integer,
but the only vectors satisfying this are generally
$\pm e_j$ and additionally for $B_4$
the vectors $\sum_{j=1}^4 \pm \frac 1 2 e_j$.
The standard $B_n$ comes from adding a simple root of the first form 
to a $D_n$ root system, 
while it is easy to check the rows of the matrix
\[
\mat{\phantom{-}1&-1&\phantom{-}0&\phantom{-}0\\
	\phantom{-}0&\phantom{-}1&-1&\phantom{-}0\\
	\phantom{-}0&\phantom{-}0&\phantom{-}1&-1\\
	-\frac 1 2&	-\frac 1 2&	-\frac 1 2&\phantom{-}\frac 1 2
}
\]
are also simple roots for a $B_4$ root system.
\epf

The union of these two lattices contains an $F_4$
root system
\[
\mat{
	0&\phantom{-}1&-1&\phantom{-}0\\
	0&\phantom{-}0&\phantom{-}1&-1\\
	0&\phantom{-}0&\phantom{-}0&\phantom{-}1\\
	\frac 1 2&	-\frac 1 2&	-\frac 1 2&	-\frac 1 2
}
\]
and so generates an $F_4$ lattice.
	Indeed, there are two distinct $\Spin(9)$ subgroups $\Kpm$
	of the group $G = F_4$ meeting in a $\Spin(8) = H$ 
	and witnessing this root data~\cite[Table~E, p.~125]{grovewilkingziller2008}.
	The resulting double mapping cylinder is $S^{25}$.

\blem\label{thm:AGbasisF4}
The family of Weyl groups $(WG,WK^-,WK^+,WH,1)$ 
corresponding to the cohomogeneity-one action
in the preceding paragraph
meets the hypotheses of \Cref{thm:AGbasis}.
\elem
\bpf
	We note that $F_4$ is simply-connected.
	The coset condition of \Cref{thm:AGbasis} is satified automatically
	if, in that notation, one of $W_j$ and $W_k$ contains the other,
	so we only need to check that for $W_j = WK^-$ and $W_k = WK^+$,
	we can take $W_\ell = WG$.
	But, as is easy to ask a computer to check~\cite{carlson2018roots}, 
	if one chooses the positive roots $\Phi^+ F_4$
	of the $F_4$ root system to be 
	$e_j$,\ \ $e_j \pm e_k$, and $\frac 1 2(1,\pm 1,\pm 1, \pm 1)$
	and the positive roots $\Phi^+ \Kpm$ 
	of the smaller groups to be subsets of these,
	then the sets
	$\{w \in WF_4 : w\Phi^+ \Kpm \subn \Phi^+ F_4\}$
	of coset representatives of $WF_4/W\Kpm$
	meet only in the neutral element.
\epf

We will need to apply \Cref{thm:AGbasis} to one other case, 
the system of subgroups of $\Sp(3)$
given by the block-diagonal subgroups 
$K^- = \Sp(2) \+ \Sp(1)$ and $K^+ = \Sp(1) \+ \Sp(2)$,
which meet in the diagonal $H = \Sp(1)^{\+3}$.
All share as a maximal torus $T = \U(1)^{\+3}$.
	It is easy to see that the roots of the larger groups in $T$
	generate an $C_3$ lattice,
	and under the standard identification of $W\Sp(3)$
	with $\Sigma_3 \semidirect \{\pm 1\}^3 < \Aut \R^3$,
	the subgroups $WK^-$ and $WK^+$
	become respectively 
	$\ang{(1\ 2)} \. \{\pm 1\}^3$ and
	$\ang{(2\ 3)} \. \{\pm 1\}^3$ ,
	while $WT$ is simply $\{\pm 1\}^3$.
	
\blem\label{thm:AGbasisC3}
The family of Weyl groups $(WG,WK^-,WK^+,WH,1)$ 
corresponding to the cohomogeneity-one action
in the preceding paragraph
meets the hypotheses of \Cref{thm:AGbasis}.
\elem
\bpf
	Note that $\Sp(3)$ is simply-connected.
	As before, the only pair of containment-incomparable
	subgroups under consideration is $\{WK^-,WK^+\}$,
	and one checks \cite{carlson2018roots}
	the sets of coset representatives
	$\{w \in WC_3 : w\Phi^+ \Kpm \subn \Phi^+ C_3\}$
	for $WC_3/W\Kpm$ meet only in $1$.
\epf

 \section{The case when one sphere 
 	is odd-dimensional}\label{sec:odd}

We now put the algebra of the previous section to use
to obtain specializations of \Cref{thm:mainK}.
In this section, at least one of the homogeneous spheres
$\Kpm/H$ is odd-dimensional.

\Kodd

\brmk
	In terms of representations,
	$\defm t$ is the class of the representation 
	$\Kp \to (\Kp)\ab/H\ab \simto \U(1)$,
	and similarly for $\defm{t_\pm}$.
	Likewise, $\defm{\ol \rho}$
	is the reduction $\rho - \dim \rho$
	of a complex $K^+$-representation
	$\rho\: K^+ \lt \U(V_\rho)$,
	trivial when restricted to $H$,
	such that the class $\defm {\ol\b}(\rho)$ 
	represented by the composition $K^+/H \to \U(V_\rho) \inc \U$
	generates $K^1(K^+/H)$,
	and similarly for $\defm{\ol\rho_\pm}$.
\ermk

\bpf[Proof of \Cref{thm:Kodd}]
  We use the description of $K^*_G (M)$ given in Theorem \ref{thm:mainK}. 
 	In both cases,  $K_G^1 (M) = 0$ 
 	since $R\Kp \lt RH$ is surjective,
 	so $
 	\KG (M) = 
 	K_G^0 (M) \iso 
 	\xu{RH}{R\Km\mn}{R\Kp}.
 	$
 	
 	\medskip
 	
 	(a)\, 
 	Recall from \Cref{thm:Keven} that $R\Km \lt RH$ is an injection and from Propositions \ref{thm:Koddcircle} and \ref{thm:Koddsphere}
 	that the map $R\Kp \simto RH[\bar\rho] \to RH$ 
 	or 
 	$R\Kp \simto RH[t^{\pm 1}] \to RH$ 
 	is reduction modulo $(\bar\rho)$ or $(t-1)$.
 	We prove the latter case; 
 	the former is similar.
 	Then the fiber product is the subring of 
 	$RH[t^{\pm1}] \x R\Kp$
 	consisting of the direct summands 
 	$\big\{(\s,\s) \in R\Kp \x R\Kp\big\}$
 	and
 	$(t-1) RH[t^{\pm 1}] \x \{0\}$.
 	We may identify the former with
	 	$R\Kp < RH < RH[t^{\pm1}]$
 	 and the latter with
	 	$(t-1) RH[t^{\pm1}] \idealneq RH[t^{\pm1}]$ 
 	and the two interact multiplicatively via the rule
 	\[
	 	\s\. (t-1)f \longbij 
	 	(\s,\s)\.\big((t-1)f,0\big) = 
	 	\big((t-1)\s f,0\big) \longbij 
	 	(t-1)\s f.
	 \]
 	%
 	
 	\medskip
 	
 	(b) 
 	We use \Cref{thm:Keven} to make identifications 
 	$R\Km \iso RH[t^{\pm1}]$ and
	$R\Kp \iso RH[\bar\rho]$
 	such that $R\Km \lt RH$ is reduction modulo 
 	$\defm{\bar t} = t-1$
 	and $R\Kp \lt RH$ modulo $\bar\rho$;
 	the other cases are the same, \emph{mutatis mutandis}.
	The fiber product can be identified as the subring of 
	$RH[t^{\pm1}] \x RH[\bar\rho]$ 
 	comprising the three direct summands
 	\[
 	\big\{(\s,\s) \in RH \x RH \big\},\qquad\qquad
 	\bar t RH[t^{\pm1}] \x \{0\},\qquad\qquad
 	\{0\} \x \bar\rho RH[\bar\rho].
 	\]
 	Multiplication across summands 
 	is determined by the three rules 
 	\[
 	(\s,\s) \. (\bar t \mspace{-1mu}f^-,0) = (\bar t\mspace{-1mu}f^-\s,0),\qquad
 	(\s,\s) \. (0,\bar\rho f^+) = (0,\bar\rho f^+\s),\qquad
 	(\bar t \mspace{-1mu}f^-,0)\.(0,\bar\rho f^+) = (0,0),
 	\]
 	so the map to $RH[t^{\pm 1},\bar\rho]/(\bar t \bar\rho)$
 	sending $(\s + \bar t\mspace{-1mu}f^-,\s + \bar\rho f^+)$
	to the class  
	 	$\s + \bar t \mspace{-1mu}f^- + \bar\rho f^+ \pmod{\bar t \bar\rho}$
 	is a ring isomorphism.
 \epf

\brmk
	This statement is obviously not the most one can say, 
	in that it can be extended
	using the extraneous description \eqref{eq:SOmap}
	of $R\SO(2n) \lt R\SO(2n-1)$
	in the proof of \Cref{thm:spheresurj} to cover the cases where
	the image of one or more of $K^\pm \lt \Homeo K^\pm/H$
	comes from an $\SO(\even)$ subgroup of $K^\pm$\,---\,%
	but this is left as an exercise 
	for the interested reader, if any,
	the current statement being long enough as it is.
\ermk

\begin{example}
	Let $M$ be the double mapping cylinder 
	associated to a diagram with 
	$H = \Spin(7)$ included in 
	$K^- = \Spin(8)$ via the standard inclusion
	and in $K^+ = \Spin(9)$ via 
	the nonstandard embedding with $K^+/H = S^{15}$;
	the larger group $G$ can be anything large enough,
	say $F_4$ or $\Spin(8) \x \Spin(9) = K^- \x K^+$.
	Then we have an explicit presentation
	\[
		\KG(M) 
			\iso 
		\Z[\s,\l^2\s,\D,\bar\rho_-,\bar\rho_+]
			/
		(\bar\rho_- \bar\rho_+)\mathrlap,
	\]
	where in $R\Spin(8) \x R\Spin(9)$, the generators
	are represented by
\eqn{
		\s 		&\longbij (\s-1,\D-\s)\mathrlap,\\
		\l^2\s	&\longbij (\l^2\s-\s-1,\l^2(\s-1)+\s-\D)\mathrlap,\\
		\D		&\longbij (\D_-,\s-1)\mathrlap,\\
	\bar\rho_-	&\longbij (\D_+-\D_-,0)\mathrlap,\\
	\bar\rho_+	&\longbij \big(0,\l^3(\s-1) - (\D-\s-1)(\s-1)\big)\mathrlap.
	}
	in the manner described in \Cref{rmk:geom-gen}.
\end{example}

%

\section{The case when both spheres are even-dimensional}\label{sec:even}
In this section we obtain the specialization of \Cref{thm:mainK}
where both the homogeneous spheres 
$K^\pm/H$ are even-dimensional.
Particularly,  $K^-$, $K^+$, and $H$ all have the same rank.
We will not have
to assume that $\pi_1 (\Kpm)$ is free abelian,
but only that the commutator subgroup $K'$
is the direct product of a simply-connected factor 
and a number of $\SO(\odd)$ factors. 
This is equivalent
to assuming $RK'$ is a polynomial ring~\cite{steinberg1975pittie}.

\begin{notation}
Occasionally
we will write $\defm T$ for a maximal torus of some 
connected, compact Lie group $\G$
and use the fact that $R\G \iso (RT)^{W\G}$ by restriction~\cite[\SS4.4]{atiyahhirzebruch},
where $\defm{W\G}$ is the Weyl group of $\G$.
Particularly, 
when $K^\pm/H$ are even-dimensional spheres,
$RH = (RT)^{WH}$ is of rank two over $R\Kpm = (RT)^{W{\Kpm}}$,
so $WH$ is an index-two subgroup of each of $W\Kpm$.
\end{notation}

We start with two similar reduction lemmas which will save
us time later.

\blem\label{thm:extrafactor}
Suppose $K^\pm,H$ are compact and connected
and there are groups $\ul{K}^\pm \leq K^\pm$
and $L,\ul H \leq H$ such that
$\Kpm = \ul K^\pm \x L$ and $H = \ul H \x L$
(we then write for short 
$(\Kpm,H) = \defm{(\ul K^\pm,\ul H) \x L}$),
and write $\ul M$ 
for the double mapping cylinder of $G/\ul H \longtoto G/\ul K^\pm$.
Then $\KG(M) \iso \KG(\ul M) \ox RL$.
\elem
\bpf
This follows from \Cref{thm:mainK} since the map
$RK^- \x RK^+ \lt RH$ then factors as
$(R \ul K^- \x R \ul K^+ \to R\ul H) \ox \id_{RL}$.
\epf

\blem\label{thm:extratorus}
Suppose $K^\pm,H$ are compact and connected
and there are groups $\ul{K}^\pm \leq K^\pm$
and $A,\ul H \leq H$ 
such that $A$ is a torus central in $\Kpm$
and meeting $\ul H$ in a finite subgroup
$F = A \inter \ul H$,
so that $H = \ul H \ox_F A$.
Then
$\Kpm = \ul K^\pm \ox_F A$ as well.
Writing $\ul M$ 
for the double mapping cylinder of $G/\ul H \longtoto G/\ul K^\pm$,
we have $\KG(M) \iso \KG(\ul M) \ox R(A/F)$.
\elem
\bpf
This follows from \Cref{thm:mainK} 
and \Cref{thm:Hodgkinplus} since the map
$RK^- \x RK^+ \lt RH$ then factors as
$(R \ul K^- \x R \ul K^+ \to R\ul H) \ox \id_{R(A/F)}$.
\epf

After application of these lemmas, it will follow from a case analysis
that most of the time we are in one of two special situations.
The easier of these two situations is when $K^-=K^+$.

\begin{restatable}{proposition}{oneimage}\label{thm:oneimage}
	Assume there exists $w$ in the identity component $N_G(H)_0$
	such that $K^+ = wK^-w^{-1}$,
	that $K^-/H=S^{2n}$ is a sphere of positive even dimension 
	and the left $K^-$-action is orientation-preserving.
	Then 
	\[
	\KG (\man)\iso RK^- \ox \K(S^{2n+1})\mathrlap.
	\]
\end{restatable}
\bpf
Note that in this case~\cite[p.~44]{grovewilkingziller2008},
$\man$ is $G$-diffeomorphic
to the double mapping cylinder of $G/H \longtoto G/K^-$,
so we may as well assume $K^+ = K^-$.
Then we may apply \Cref{thm:mainK},
noting that $RK^- \inter RK^+ = RK$
and that by \Cref{thm:Kevensphere},
\[
	\frac{RH}{RK^- + RK^+} = \frac{RK^-\{1,\rho\}}{RK^-} \iso RK^-\{\rho\}.
	\qedhere
\]
\epf

\brmk
Forgetting the manifold itself
and proceeding in terms of representation theory,
we could also have noted
that if $K > H$ share a maximal torus and $w$ lies in $N_G(H)_0$,
then $wKw^{-1}$ also contains that torus,
with respect to which $WK = W(wKw^{-1})$.

Proceeding more topologically, on the other hand, we could note that
if $K^+ = K^- = K$, 
then the natural map $BH \lt BK$
allows us to define a sphere bundle
$S(K/H) \to M_G \to BK$. 
The proof of the analogue for Borel cohomology~\cite[Prop.~5.2]{CGHM2018}
worked by showing this bundle was cohomologically trivial,
and it is to reflect this analogy that we retain the number $n$.
\ermk

\brmk
It is interesting to note that if we do not have $K^+ = K^-$, 
then $H = K^- \inter K^+$.
To see this, first note that 
since $K^- \inter K^+$ and $H$ share a maximal torus,
$(K^- \inter K^+)/H$ is even-dimensional.
But $(K^- \inter K^+)/H \to K^+ / H \to K^+ /(K^- \inter K^+)$
is a fibering of a sphere over a simplicial complex
and by connected simplicial complexes,
and Browder showed that when the fiber is none of $S^1$, $S^3$, or $S^7$, 
either the base or the fiber of such a bundle 
must be trivial~\cite{browder1963higher}.

But this dichotomy does not 
lead to a dichotomy in expressions for $\KG(M)$.
For example, the block-diagonal subgroup $H = \SO(4) \+ [1]^{\oplus 2}$
of $G = \SO(6)$ is the intersection of
$K^- = \SO(5) \+ [1]$ and $K^+\mn = w K^- w^{-1}$
for $w =[1]^{\oplus 4} \oplus 
\big[\begin{smallmatrix} 0 & 	\!\!\!\phantom{0}-1\\1 & 
								\!\!\!\phantom{-1}0\end{smallmatrix}\big]$,
which lies in $[1]^{\oplus 4} \+ \SO(2) < N_G(H)_0$.
Thus, up to diffeomorphism, 
the inclusion diagram $(G,K^-,K^+,H)$ 
expresses the same double mapping cylinder $M$
as the one instead taking $K^+ = K^- = \SO(5)\+[1]$.
\ermk

The other easy-to-manage special case follows from 
a less trivial product decomposition.

\bprop\label{thm:twofactor}
Let connected, compact Lie groups $\ul K^\pm > \ul H^\pm$ 
be such that $\ul K^\pm/\ul H^\pm = S^{2\npm}$ 
are even-dimensional spheres.
Write $H = \ul H^- \x \ul H^+$ and consider it in the natural
way as a subgroup of 
$K^- =  \ul K^- \x \ul H^+$, of
$K^+ = \ul H^- \x \ul K^+$, and of
$G = \ul K^- \x \ul K^+$.
Then if $M$ is the double mapping cylinder of $G/H \longtoto G/\Kpm$,
we have 
\[
	\KG M \iso RG \ox \ext[z] 
\]
for a generator $z$ of degree $1$.
\eprop
\bpf
By \Cref{thm:Kevensphere}, 
we know $R\ul H^\pm$ is free of rank two over $R\ul K^\pm$,
say on bases $\{1,\s_\pm\}$.
Then $RK^-$, $RK^+$, and $RH$ are free over $RG = R\ul K^- \ox R\ul K^+$
respectively on the bases 
\[\{1\ox 1,\,\s_- \ox 1 \}, \qquad 
 \{1\ox 1,\,1 \ox \s_+\}, \qquad 
  \{1\ox 1,\ \s_- \ox 1,\ 1 \ox \s_+,\ \s_- \ox \s_+\}\mathrlap.\]
Thus, by \Cref{thm:mainK},
we see $K^0_G (M)$ is the intersection of $R\Kpm|_H$,
namely the free $RG$-module on $1 \ox 1$,
and $K^1_G (M) \iso RH/(RK^- + RK^+)$ 
is the free cyclic $RG$-module on $z = \d(\s_- \ox \s_+)$.
Thus $\KG(M)$ is a free $RG$-module on $1 \in K^0_G(M)$ and $z \in K^1_G(M)$,
and since $2z^2 = 0$ by antisymmetry and $\KG(M)$ is torsion-free,
it follows $z^2 = 0$.
\epf

\brmk
The manifold $M$ is a sphere $S^{2n_- + 2n_+ + 1}$ 
under these conditions.\footnote{\ This will also hold
	if either sphere or both is odd-dimensional.} 
Indeed, the fiber over $-1$ is $S^{2n_-}$, 
that over $1$ is $S^{2n_+}$, 
and in the interior the fiber is the product
of the two, so $M$ is the join 
$S^{2n_-} * S^{2n_+}$.
\ermk

\bex[{\cite[Sect.~4.3]{puettmann2009}}] We use \Cref{thm:oneimage} to compute the equivariant cohomology of the space $M$ arising from the inclusion diagram 
\[(G,\Km,\Kp,H) = 
	\big(
		\Sp(2), 
		\Sp(1)^2, 
		\Sp(1)^2,
		\Sp(1) \x \U(1)
	\big)
	\mathrlap.
\]
P\"uttmann shows
$\H(M;\Z) \iso \H(S^3;\Z) \ox \H(S^4;\Z)$
using the \MVS,
so from the \AHSS we see $\K(M) \iso \K (S^3) \ox \K (S^4)$ as well.
The restriction of the defining representation $\s$ of 
$\Sp(1) < \Quat^\x$ on $\Quat \iso \C \+ j\C$ 
to the maximal torus $\U(1) < \C^\x$ is $t + t\-$,
where $t$ is the defining representation,
so 
\[
	K^1_G (M) \iso
	\frac{\Z[\s] \ox \Z[t^{\pm 1}]}{\Z[\s] \ox \Z[t+t\-]}
		\iso
	\Z[\s] \ox t \Z[t+t\-]
	\iso R\big(\Sp(1)^2\big)[1]
\]
as expected. 

This action is equivariantly formal for Borel cohomology with 
integer coefficients~\cite[Cor.~1.3]{goertschesmare2014},
and from \Cref{thm:eqf}, 
it is equivariantly formal for $\KG$ too,
but it is illuminating to show this explicitly
by examining the forgetful map $\KG \lt K$ 
on the \MVS of the standard cover.
By the snake lemma, 
this amounts to checking the maps
\[
R\G \isoto K^0_G(G/\G) \lt K^0(G/\G)
\]
taking a representation $V_\rho$ of $\G$
to the bundle $G \ox_\G V_\rho \lt G/\G$
are surjective for $\G \in \{K^\pm,H\}$.\footnote{\ 
In fact, applying the module structure in 
\Cref{thm:mainK}
to both sequences, it would be enough just to see
$K^0_G M \lt K^0 M$
is surjective, 
and once we know $K^1(G/H) = K^1 \CP^3 = 0$,
it would suffice to prove $RK \lt K^0(G/K)$
is surjective,
but the same proof involves both maps.}
It is not hard to check this map takes 
$1\ox t \in R\big(\Sp(1) \x U(1)\big)$
to the tautological bundle $\g$ over $\CP^3$
and $1 \ox \s \in R\big(\Sp(1)^2\big)$
to the tautological bundle $\xi$ over $\Quat P^1$.\footnote{
Note $\Sp(2) \lt S^7$ given by 
$\smash{A \lmt A\.\bigl[\begin{smallmatrix}
0\\1 
\end{smallmatrix}\bigr]}$
has stabilizer $\Sp(1) \+ 1$
and transforms the action of
$1 \+ \Sp(1)$ to scalar right-multiplication on $S^7 \subn \Quat^2$,
so the total spaces of the bundles are $S^7 \ox_{\Sp(1)} \Quat$
and $S^7 \ox_{\U(1)} \C$.}
Since $\H(\CP^3) = \Z[c]/(c^3)$,
where $c = c_1(\g)$,
and $c_1$ induces an isomorphism $\wt K^0 (\CP^1) \isoto H^2 (\CP^1)$,
this gives us surjectivity for $H$.
As for $K^\pm$,
since $\s$ restricts to $\U(1)$ as $t + t\-$,
we see the pullback of $\xi$ over $\CP^3$ is $\g \+ \g^\vee$.
%
%
The total Chern class $1 + c_2(\tau) \in \H(\Quat \mr P^1)$
hence pulls back to $(1+c)(1-c) 
								\in \H (\CP^3)$.
The \SSS of $S^3/S^1 \to \CP^3 \to \Quat \mr P^1$
collapses for degree reasons,
so that $H^4 (\Quat P^1) \to H^4 (\CP^3)$.
Thus, since $-c^2$ generates $H^4 (\CP^3)$, 
also $c_2(\tau)$ generates $H^4 (\Quat P^1)$.
As \[\wt K^0 (S^4) \iso \wt K^4 (S^4) \iso \wt K^0 (S^0) = \Z\] 
and the Chern character induces a natural isomorphism $\K \ox \Q \lt \H(-;\Q)$
on finite complexes, it follows 
$[\tau]$ generates $\wt K^0(S^4)$ as needed.
\eex

The desired simultaneous generalization of
\Cref{thm:oneimage,thm:twofactor} is as follows.

\Keven

The proof 
has been factored into as many 
Lie-theoretic lemmas and reduction steps
as possible but still seems
to unavoidably be a bit of a slog.

\bpf[Proof of \Cref{thm:Keven}]
Note that the images $\defm\Kepm$ of the action maps
$\defm{\apm}\: \Kpm \lt \Homeo \Kpm/H$
are by definition effective and hence must be  
$\SO(2n+1)$ or $G_2$, 
with the image of $H$ being $\SO(2n)$ or $\SU(3)$ respectively%
~\cite[Ex.~7.13]{Be}\cite[Table~C, p.~104]{grovewilkingziller2008}.
The effective images $\defm{\Hepm} = \apm(H)$ of $H$, 
in particular, 
determine $\Kepm$ uniquely up to isomorphism.

Most of the proof involves analyzing the configurations
of these preimages after stripping away
extra tensor factors
to eventually arrive at a base case.
The recurrent phrase ``\defd{factor out $\Pi$}'' means
to apply \Cref{thm:extrafactor}
and analyze the remaining system of isotropy groups 
$\ul K^- \from \ul H \to \ul K^+$,
whereas ``\defd{factor out $A/F$}'' means to apply \Cref{thm:extratorus}.
We say we have reduced to a \defd{\casep}
if \Cref{thm:twofactor} applies, 
in which case that branch of the case analysis terminates,
and similarly say we have reduced to a \defd{\caseb}
if \Cref{thm:oneimage} applies.
Beyond these base case schemata, 
there are a few exceptional base cases 
enumerated in \Cref{sec:even-lemma},
which as we have mentioned, 
all turn up as examples in the literature,
and the case with
$\He^- \iso \SO(2) \iso \He^+$.

\bigskip

\nd\textbf{0.}\ \emph{The case neither of $\Hepm$ is a circle}

\medskip

As $K^\pm/H$ are even-dimensional spheres of dimension $ > 2$,
the long exact fibration sequence of $H \to \Kpm \to \Kpm/H$
induces isomorphisms $\pi_1 H \isoto \pi_1 \Kpm$.
It follows that the inclusion of $A = Z(H)_0$ in $H$
induces surjections $\pi_1 A \lt \pi_1 \Kpm$
and we can write $\Kpm$ as $(\Kpm)' \ox_{(F^\pm)} A$
for $\defm{F^\pm} = \ker\big((\Kpm)' \x A \epi \Kpm\big)$.
Since $\Kpm/H$ are spheres, 
by two applications of \Cref{thm:lowerF} we have
$H' \ox_{F^-} A = H = H' \ox_{F^+} A$,
so $\defm F = F^- = F^+$.
Thus the inclusions $H \longtoto \Kpm$
factor as virtual product maps of the form $i_{\pm} \ox_F \id_A$.
Factoring out $A/F$,
we need only analyze 
$\KG(M')$
for $M'$ the double mapping cylinder of $G/(\Kpm)' \longtoto G/H'$.
We may thus adopt the notational convenience of 
assuming the groups $\Kpm$ of the original
triple $(\Kpm,H)$ were semisimple.

Recall that a closed subgroup of a simply-connected Lie group
can be written as a direct product of closed subgroups 
of its simple factors~\cite[p.~205]{boreldesiebenthal}
and that normal subgroups can be written as products 
of simple factors and finite central groups.
Examining $\apm$ and $\apm|_H$
on the Lie algebra level,
we see their kernels contain all but one
of these simple factors, 
or all but two in case
$\Hepm = \SO(4) \iso \SO(3)^2 / \big\{\mn\mn\pm\mn(I,I)\big\}$
is not simple.
Thus we have product decompositions
\eqn{
	\Kpm 	&\iso \tKepm \x \Pi^\pm\mathrlap,\\
	H 		&\iso \tHepm \x \Pi^\pm\mathrlap,
}
where 
the ineffective kernels
$\defm{\Pi^\pm} \ceq \ker \apm$ are
products of simply-connected and $\SO(\odd)$ factors,
their complements
$\defm{\tKepm} \leq \Kpm$
induce isomorphisms or double-coverings $\tKepm \inc \Kpm \epi \Kepm$,
and $\defm{\tHepm}$ are the intersections of $H$ and $\tKepm$,
accordingly singly or doubly covering $\Hepm$ under $\apm$.


\smallskip


\bitem
	\item \emph{Suppose it is possible to select $\tKepm$
		in such a way that 
		$\tHe^+ = \tHe^- \eqc \defm{\tHe}$.}

	\smallskip

	Then $\Pi^+ = \Pi^-$
	and we may factor it out.
		What remains is the pair of inclusions $\tHe \longtoto \tKepm$, 
		so we examine the images of $R\tKepm \longtoto R\tHe$.
	\bitem
		\item \emph{Suppose that $\tHe \niso \Spin(8)$.}
		
		\smallskip
		
		An inclusion $\SO(2n) \longinc \SO(2n+1)$ for $n \neq 4$
		or $\SU(3) \longinc G_2$
		induces an inclusion of root lattices in a unique way
		by \Cref{thm:latticeinclusionG2,thm:latticeinclusionBn}.
		It follows that the maps $R\tKepm \longtoto R\tHe$
		have the same image, so we have a \caseb.
	
		\item \emph{Suppose that $\tHe \iso \Spin(8)$.}
	
	\smallskip
		\bitem
			\item
			If the inclusions of root lattices
			induced by $\tHe \longtoto \tKepm$ 
			are both standard,
			then as in the previous item, 
			we have a \caseb.

		\smallskip
	
			\item

		\smallskip

			Otherwise our $B_4$ lattices are both of 
			those described
			in \Cref{thm:latticeinclusionBn}
			and so together span an $F_4$ lattice,
			and the intersection 
			$RK^- \inter RK^+$ in $RH = R\Spin(8)$ is $RF_4$.
			By \Cref{thm:AGbasisF4}, then,
			$R\Spin(8)$ is free over $RF_4$
			on $1152/192 = 6$ elements
			and each $R\Spin(9)$ 
			is free on $1152/384 = 3$ elements,
			so by arithmetic,
			\[
			\frac{R\tHe}{R\tKe^- + R\tKe^+}
				\iso 
			RF_4 \iso R\tKe^- \cap R\tKe^+
				\mathrlap.
			\]
		\eitem
	\eitem
	\item \emph{Suppose it is impossible to select $\tKepm$
		in such a way that $\tHe^- = \tHe^+$.}
	\bitem
		\item
		\emph{Suppose that neither of $\Hepm$ is isomorphic to $\SO(4)$.}
	
		\smallskip
	
	The assumption implies $\Hepm$ 
	and hence the single or double covers $\tHepm$ are simple.
	Since $H$ is a product of simply-connected groups and $\SO(\odd)$
	factors, and since subgroups $\tKepm \leq \Kpm$ 
	singly or doubly covering $\Kepm$ under $\apm$
	cannot be chosen such that $\tHepm = H \inter \tKepm$
	agree, we must have $\tHe^- \inter \tHe^+ = 1$.
	Thus there exists a factorization
	\[
		H = \tHe^- \x \tHe^+ \x \Pi
	\] 
	for $\Pi$ a product of totally ineffective factors
	contained in $\Km \inter \Kp$.
	Since $\rk \tKepm = \rk \tHepm$ 
	and the groups $\tHepm$ are simple, it follows 
	\[
	\tKe^+ \inter \tHe^- = 1 = \tKe^- \inter \tHe^+\mathrlap,
	\]
	and as $H = \tHe^- \x \tHe^+ \x \Pi$ is contained in 
	both groups
	$\Kpm$, 
	they must admit abstract decompositions
	\eqn{K^- &\iso 
				\tKe^- \x \tHe^+ \x \Pi\mathrlap,\\
		K^+ &\iso	
				\tHe^- \x \tKe^+ \x \Pi
	}
	respecting the inclusions $\tHepm \longinc \tKepm$.
	Thus we may factor out $R\Pi$ and 
	afterwards have a \casep. 
	
	\smallskip

	\item \emph{Suppose
		at least one of $\Hepm$ is isomorphic to $\SO(4)$.}
	
	\smallskip
	
	We may suppose without loss of generality that it is
	$\smash{\He^+}$ which is isomorphic to $\SO(4)$, 
	so that ${\tHe^+} \iso \Spin(4) \iso \Sp(1)^2$
	and $\tKe^+ \iso \Spin(5) \iso \Sp(2)$.
	Since $\tHe^+$ and $\tHe^-$ are both direct factors
	of the semisimple group $H$
	and we have assumed that $\tHe^- \neq \tHe^+$,
	we have a dichotomy based on whether $\tHe^-$
	shares $0$ or $1$ of the $\Sp(1)$ factors of $\tHe^+$.
	\bitem
	\item
	\emph{Suppose no $\Sp(1)$ factor of $\tHe^+$ lies in $\tHe^-$.}
	
	\smallskip
	
	Then $\tHe^+ \leq \Pi^-$, so we have 
	\[
	H \iso \tHe^- \x \Pi^- \iso \tHe^- \x \tHe^+ \x L 
	\]
	for some direct complement $L$ with $\Pi^- \iso \tHe^+ \x L$.
	It follows
	\[
	K^- \iso \tKe^- \x \tHe^+ \x L\mathrlap.
	\]
	On the other hand,
	the inclusion $H \longinc K^+$
	factors abstractly as 
	\[\tHe^- \x \tHe^+ \x L \longmono \tKe^+ \x \Pi^+\mathrlap,\]
	with the image of $\tHe^+$ lying in $\tKe^+$,
	so it follows $\Pi^+ \iso \tHe^- \x L$.
	Thus we factor out $L$ and achieve a \casep.

	\smallskip
	
	\item
	\emph{Suppose one $\Sp(1)$ factor of $\tHe^+$ lies in $\tHe^-$.}

\smallskip
	Since $\He^-$ is isomorphic to either $\SU(3)$ or $\SO(\even)$
	and $\tHe^-$ is a product of direct factors of $H = \Sp(1)^2 \x \Pi^+$,
	we must also have $\tHe^- \iso \Sp(1)^2$
	and $\tKe^- \iso \Sp(2)$.
	Factoring out $\Pi^- \inter \Pi^+ < H$, 
	what remains are the inclusions $\tH \longtoto \tKepm$,
	which can be identified with 
		\[
			\Sp(2) \x \Sp(1) \longfrom \Sp(1)^3 \longto \Sp(1) \x \Sp(2).
		\]
		Then by \Cref{thm:AGbasisC3}, 
		$R\Sp(3)$ is free over $R\big(\Sp(1)^3\big)$
		on $6 = |\Sigma_3|$ elements
		and each of $R\tKe^\pm$ 
		is free on $3$ elements,
		meaning
		\[
	\frac{R\tHe}{R\tKe^- + R\tKe^+} \iso RC_3 \iso R\tKe^- \inter R\tKe^+
		\] 
		as expected.
	\eitem
	\eitem
%
		
\eitem

\nd\textbf{1.}\ \emph{The case exactly one of $\Hepm$ is a circle}
	
	\medskip

Without loss of generality, assume that $\He^- \iso \SO(2)$
and $\He^+ \niso \SO(2)$.
As before let $\defm\tKepm$ be complements to the normal subgroups
$\ker \apm \lhd \Kpm$ and $\defm\tHepm = H \inter \tKepm$.
By our assumption on the structure of $K^-$,
we can write
\[
	K^-	\iso \xt F {(\tKe^- \x \Pi^-)}  A 
\]
for $\defm A = Z(K^-)_0$
and $\defm{\Pi^-}$ a direct complement to $\tKe^-$ in 
the commutator group $(K^-)'$,
and $\defm F \iso (\tKe^- \x \Pi^-) \inter A$.
Then $H \inter (\tKe^- \x \Pi^-) = \tHe^- \x \Pi^-$,
and $K^-/H \homeo S^2$,
so by \Cref{thm:lowerF}, 
we may write $H \iso (\tHe^- \x \Pi^-) \ox_F A$.
Since $\tHe^-$ is a circle, we have $H' = \Pi^-$.

Now $\tKe^+$ is not isomorphic to either $\Spin(3)$ or $\SO(3)$,
so $\tHe^+$ is a closed subgroup of $\Pi^-$.
By our assumption on $(K^+)'$, then,
$\tKe^+$ is a direct factor and there exists a complement 
$\defm L \lhd \Pi^-$
with 
\eqn{
	\Pi^- &\iso L \x \tHe^+\mathrlap,\\
	(K^+)' &\iso L \x \tKe^+\mathrlap.
	}
It is clear then that $K^+ = \tHe^-\.(L \x\tKe^+)\. A$.
We have
\[
\tHe^- \inter (L \x \tKe^+) = 
\tHe^- \inter H \inter (L \x \tKe^+) =
\tHe^- \inter (L \x \tHe^+) = 1
\]
and also
\[
(\tHe^- \x L \x \tKe^+) \inter A = 
(\tHe^- \x L \x \tKe^+) \inter H \inter A = 
(\tHe^- \x L \x \tHe^+) \inter A = 
F\mathrlap,
\]
so in fact $K^+ \iso (\tHe^- \x L \x \tKe^+) \ox_F A$.
%
%

Thus we may factor out $A/F$ and then $L$ to obtain a \casep.

\bigskip
		
\nd\textbf{2.}\ \emph{The case $\Hepm$ are both circles}
		
		\medskip

The intersections  $\defm{\Pi^\pm}$ of 
$(\Kpm)'$ with the ineffective $\ker \apm$
admit complements $\tKepm$ in $(\Kpm)'$ by assumption.
Since $\im \apm \iso \SO(3)$ is simple and centerless,
the centers $Z(K^\pm)$ are also contained in $\ker \apm$.
This kernel is obviously contained in the stabilizer $H$
as well,
so $\Pi^\pm = (\Pi^\pm)' \leq H'$.
On the other hand, 
since the images $\apm(H) \iso \SO(2)$ are abelian,
the commutator subgroup $H'$ is contained in 
both of $\ker \apm$,
so $\Pi^\pm = H'$.

By the assumption on $(\Kpm)'$,
we have
\eqn{
	\Kpm 	&\iso (H' \x \tKe^\pm) \. Z(\Kpm)_0\mathrlap,\\
	H 		&\iso (H' \x \smash{\tHe^\pm}) \. Z(\Kpm)_0\mathrlap.
}
Now 
consider
the torus $\defm A \ceq \big(Z(K^-) \inter Z(K^+)\big)_0$.
Taking $\defm{\ul H} = H'  \tHe^- \tHe^+$,
and $F = \ul H \inter A$,
we may write $H \iso \ul H \ox_F A$.
If we set $\defm{\ul K^\pm} = \ul H \tKepm$,
then evidently $\ul K^\pm \inter H = \ul H$ and $\Kpm = \ul K^\pm  A$.
Since 
\[
\ul K^\pm \inter A = \ul K^\pm \inter H \inter A = \ul H \inter A = F\mathrlap,
\]
we find $\Kpm \iso \ul K^\pm \ox_F A$,
so we may factor out $A/F$.

	\bitem
		\item 
	\emph{Suppose $A = Z(K^-)_0 = Z(K^+)_0$.}

	\smallskip
	In this case $Z(H)/A$ is one-dimensional,
	so we may select $\tKepm$
	in such a way that $\defm \tHe = \tHe^- = \tHe^+ \iso \SO(2)$.
	Factoring out $A/F$ and then $H'$
	leaves a configuration
	$\SO(2) \longtoto \tKepm$
	where $\tKepm$ are each $\SO(3)$ or $\Spin(3)$.
	Either way, the induced
	map $R\tKepm \lt R\SO(2) = \Z[t]$
	has image $\Z[t+t\-]$,
	so we are functionally in the situation of \Cref{thm:oneimage}
	and in particular 
	\[
	\frac{R\tHe}{R\tKe^- + R\tHe^+} \iso 
	\frac{\Z[t]}{\Z[t+t\-]}			\iso
	t\.\Z[t+t\-]
	\]
	is of rank one over $\Z[t+t\-]$.

	\smallskip
	
		\item
	\emph{Suppose $Z(K^-)_0 \neq Z(K^+)_0$.}
	
	\smallskip 
	
	Write $\defm{T}$ for the two-dimensional torus
	$\tHe^-\.\tHe^+$ in $H$.
	Then after factoring out $A/F$ we have to deal with
	the inclusions of $H = H' \x T$ 
	in $(H' \x \tKe^\pm) \. S^1$,
	where $\id_{H'}$ factors out of these inclusions
	but we claim nothing particular about 
	the two inclusions $T \longmono \tKe^\pm\. S^1$.
	Factoring out $H'$,
	we arrive at $\ul H = T$ and 
	$\ul K^\pm \iso \tKepm \ox_{F} S^1$,
	where $|F| \leq 2$.
	
	The inclusions $T \longinc \ul K^\pm$
	induce inclusions $R\ul K^\pm \iso (RT)^{\ang{\wpm}} \longinc RT$,
	where $\wpm$ generates $W\ul K^\pm \iso \Z/2$.
	Identifying $RT^2$ with the group ring $\Z\X$
	of the character group 
	$\defm\X = \defm{\X(T)} = \Hom(T,S^1)$,
	these can be seen as induced by two reflections
	of the vector space $\ft\dual \iso \R^2$ 
	which preserve the integer lattice $\X(T) \iso \Z^2$.
	Under this identification $\defm W = \ang{w_-,w_+}$ becomes a dihedral
	subgroup $D_{2k}$ of $\GL(2,\Z)$.
	These are classified: they can only be $D_4, D_6, D_8, D_{12}$
	and are conjugate to the standard presentations 
	for the Weyl groups of types
	$D_2 = A_1 \x A_1$,\  $A_2$,\  $B_2 = C_2$,\  and $G_2$
	as well as a second $D_6 < WG_2$
	not generated by root reflections,
		which hence does not occur~\cite[Prop.~1]{tahara1971finite}\cite{mackiw1996finite}.
	The root lattice $Q_W$ and weight lattice $P_W$ 
	corresponding to reflection groups $W$
	of this type in $\R^2$ are unique 
	(up to equivariant isomorphism) and there are examples,
	most of which we produce immediately 
	following the present argument,
	showing any intermediate lattice between $Q_W$ and $P_W$
	occurs as $\X$ for some cohomogeneity-one action.
	
	In all of these cases, we need to see
	\[
		\defm \Theta 
			\ceq 
		\frac{RT} 
			{(RT)^{\ang{w_-}} + (RT)^{\ang{w_+}}}
	\]
	is a free cyclic module over $(RT)^W$.
	One is tempted is to use \Cref{thm:AGbasis},
	but it can happen that $RT$ is not free over $(RT)^W$.
	Instead our answer comes from the Stiefel diagram.
	The ring $RT$ is free on the $\Z$-basis $\X$.
	Quotienting by $(RT)^{\ang{w_-}}+(RT)^{\ang{w_+}}$,
	annihilates $\X^{\ang{w_-}}$ and  $\X^{\ang{w_-}}$
	and induces relations
	\eqn{
		\phantom{\qquad\qquad\qquad\mbox{for each }	
			\t \notin \X^{\ang{w_-}}\mathrlap,}
		\mathllap{	w_-} \t	 &\equiv -\t
						 \qquad\qquad\qquad\mbox{for }	
							\t \notin \X^{\ang{w_-}},\\
		\phantom{\qquad\qquad\qquad\mbox{for each }	
			\t \notin \X^{\ang{w_-}}\mathrlap,
				}
			\mathllap{	w_+} \t	 &\equiv -\t
					\qquad\qquad\qquad\mbox{for }	
							\t \notin \X^{\ang{w_+}}\mathrlap,
		}
	since $\t + w_- \t \in (RT)^{\ang{w_-}}$
	and $\t + w_+ \t \in (RT)^{\ang{w_+}}$.
	It follows $\Theta$ admits a $\Z$-basis given
	by those characters of $T$ lying in the interior $\defm C$
	of a fundamental domain.%
	\footnote{\ 
		The notation $C$ is meant to suggest a Weyl chamber,
		even though our dihedral group is just a group of symmetries of a
		$\Z^2$ lattice, not \emph{a priori} the Weyl group of anything,
		because the same reasoning goes through.
		}
	
	On the other hand, 
	$(RT)^W$ is spanned by orbit sums $\defm{S\t} = \sum_{w \in W/\Stab \t} w\t$.
	These are indexed by $W$-orbits of $\X$,
	of which there is precisely one 
	per character $\t$ in the closed fundamental domain $\defm{\ol C}$. 
	Drawing out the diagrams, one checks for each lattice type
	that there is a minimal strongly dominant integral weight $\l_0$,
	which makes $\t \longbij \t \. \l_0$
	a bijection $\ol C \inter \X \longbij C \inter \X$.\footnote{\ 
	If $\X$ is the lattice spanned by the fundamental weights
	dual to the simple roots of the root system for $W$,
	so that half the sum of positive roots is an integral weight $\rho$,
	then~\cite[Lem.~5.58]{adamsLiebook}
	we have $\rho = \l_0$.
	But these are not all the cases.}
	Recall 
	that if $\X$ is given the partial order
	determined by setting $\s \geq \t$
	just when $\t$ lies in the convex hull of 
	the orbit $W\.\s$,
	then given $\s,\t \in \X \inter \ol C$,
	the difference $S(\s\t) - S\s\.S\t$ 
	is a sum of terms of lower order~\cite[Prop.~6.36]{adamsLiebook}.
	If we filter $\Theta$ with respect to this order,
	then it follows the $(RT)^W$-module structure
	on the associated graded module $\gr \Theta$
	is given by 
	$S\s \. \ol{\t\l_0} = \ol{(\s\t)\l_0}$,
	so $\Theta$ is the free cyclic $(RT)^W$-module
	generated by $\l_0$ as claimed.\qedhere
	\eitem
\epf

\brmk\label{rmk:existingexamples}
	It is interesting to note that all of the
	exceptional cases occur as the
	``degree-generating actions'' tabulated by 
	P{\"u}ttmann~\cite[{\SS}5.2]{puettmann2009}%
	\cite[Table~E, p.~105]{grovewilkingziller2008}.
	The actions of $F_4$ on $S^{25}$
	and $\Spin(3)$ on $S^{13}$
	already came up in the ``no circular isotropy''
	case, and the others are among the ``two circles'' cases,
	as per the following examples.
\ermk

\bex
The dihedral group $D_4$, a Coxeter group of Killing--Cartan type $D_2$,
is realized as the Weyl group of a cohomogeneity-one action with $H \iso T^2$ as follows.
One has an isomorphism 
$\SO(4) \iso \xt{\Z/2}{\Spin(3)}{\Spin(3)}$
and can consider the diagram
\[
G \iso \SO(4),\quad 
K^- = \xt{\Z/2}{\Spin(3)}{\Spin(2)},\quad
K^+ = \xt{\Z/2}{\Spin(2)}{\Spin(3)},\quad
H = \xt{\Z/2}{\Spin(2)}{\Spin(2)} = T.
\]
Write $\wt T = \Spin(2) \x \Spin(2)$
and $R\wt T = \Z[s,t, s\-t\-]$.
Then $W = W\SO(4) \iso S_2 \x \{\pm 1\}$. 
Since $\SO(4)$ is not simply-connected~\cite{steinberg1975pittie},
we see $RT = \Z[s^{\pm 1}t^{\pm 1}]$ is not free over
\[
	R\SO(4)\iso 
	(RT)^W \iso
	\Z[s+s\- + t + t\-,st + s\-t\-,s\-t+st\-]
	\mathrlap,
\]
illustrating the proof of the $\Hepm = \SO(2)$ case in
\Cref{thm:Keven} cannot be run through
\Cref{thm:AGbasis} in all cases.

Instead considering the two-fold 
covers inside $G = \Spin(4) \iso \Spin(3)^2$,
one obtains a Weyl group of type $D_2$ again,
but now $R\wt T =  \Z[s,t, s\-t\-]$
is free over 
\[
	R\Spin(4) = (RT)^W \iso \Z[s+s\- + t + t\-,st + s\-t\-]
	\mathrlap,
\]
and one can apply \Cref{thm:AGbasis} again.
The space acted on is $S^2 * S^2 \homeo S^5$.

We leave it to the reader to construct an analogous example 
with $G = \SO(3) \x \SO(3)$.
\eex

\bex
The dihedral group $D_6$, a Coxeter group of Killing--Cartan type $A_2$,
is realized  as the Weyl group of a cohomogeneity-one action with $H \iso T^2$ as follows.
Consider the diagram 
\[
	G = \U(3),			\qquad\qquad\quad
	K^- = \U(2) \x \U(1),\qquad\qquad\quad
	K^+ = \U(1) \x \U(2),\qquad\qquad\quad
	H = \U(1)^3\mathrlap.
\]
	In the notation of the proof of \Cref{thm:Keven},
	the irrelevant torus 
	$A = Z\big(\U(3)\big) \iso S^1$ is the group of diagonal matrices
	and $F \iso \ang{e^{2\pi i/3}}$.
	After factoring out $A/F$,
	one has the corresponding subgroups of $\SU(3)$,
	and the manifold is $S^7$.
	The reduced $\ul\Kpm$ are both isomorphic to $\U(2)$,
	and one has $W = W\SU(3) = \Sigma_3$
	with $w_- = (1\ 2)$ and $w_+ = (2\ 3)$.
	Since $\SU(3)$ is simply-connected
	and it is easy to check the coset condition applies,
	one could also apply \Cref{thm:AGbasis}.
\eex

\bex
The dihedral group $D_8$, a Coxeter group of Killing--Cartan type $BC_2$,
is realized as the Weyl group of a cohomogeneity-one action with $H \iso T^2$ as follows.
Consider the diagram 
	\[
	G \iso \SO(5),\qquad 
	K^- = \U(2) \x \{1\},\qquad
	K^+ = \SO(2) \x \SO(3),\qquad
	H = \SO(2) \x \SO(2) \x \{1\} = T,
	\]
where all subgroups are block-diagonal,
$\U(2) \+ \{1\}$ being embedded in the block-diagonal
$\SO(4) \+ \{1\}$ in the expected manner.
	Then $WG \iso \Sigma_2 \semidirect \{\pm 1\}^2$
	is a Coxeter group of type $B_2$ 
	acting on $\ft^2$ as the dihedral group $\smash{D_8}$
	and is generated by
	$\smash{w_- = \big((1 \ 2), 1,1\big)}$
	and $\smash{w_+} = (\id,1,-1)$.
	\Cref{thm:AGbasis} does not apply
	as stated, as $\SO(5)$ is not simply-connected,
	but the relevant part of Steinberg's proof~\cite{steinberg1975pittie}
	only requires that $R\SO(5)$ be polynomial, which it is,
	and one can check the coset condition holds.

	One can also consider the cover 
	\[
		G = \Spin(5) = \Sp(2),\qquad\quad
		K^- = \U(2),\qquad\quad
		K^+ = \U(1) \+ \Sp(1),\qquad\quad
		H = \U(1) \+ \U(1) = T\mathrlap,
	\]
	which generates the same $W$.
\eex

\bex
The dihedral group $D_{12}$, 
a Coxeter group of Killing--Cartan type $G_2$,
is realized as the Weyl group of a cohomogeneity-one action with $H \iso T^2$ as follows.
Consider the adjoint action of the compact exceptional group $G_2$
on its Lie algebra $\fg_2 \iso \R^{14}$.
This restricts to an action on the unit sphere $S^{13}$
under the norm induced by the Killing form,
and the orbits are given by the intersection of $S^{13}$
with a Weyl chamber in the Lie algebra $\ft^2$ of a maximal torus,
cutting out an arc of the unit circle $S^1 \subn \ft^2$
of angle $\pi/6$.
The principal isotropy group fixing a point on the interior of the arc
is $T^2$ itself
and the singular isotropies fixing the endpoints 
are two nonconjugate copies of $\U(2)$~\cite{miyasaka2001adjointG2}.
The reflections $w_\pm$ generate the dihedral group $WG_2 = D_{12}$. 
As $G_2$ is simply-connected, one can 
check the coset condition and
apply \Cref{thm:AGbasis} again.
\eex

\section{Equivariant formality}\label{sec:eqf}

In this final section, we let $G \act M$ be a cohomogeneity-one
action with $M/G$ a closed interval 
as in the first fork \ref{fig:double} of
Mostert's dichotomy \ref{thm:Mostert}
and use the structure theorems for $\KG(M)$
in the previous two sections and the representation theory of 
\Cref{sec:mapRG} to characterize equivariant formality of such actions.

Recall that
\defd{K-theoretic equivariant formality}
means surjectivity of the map
$\smash{\KG(M) \lt \K(M)}$
forgetting the $G$-equivariant structure on a complex vector bundle.
This condition, 
first studied by Matsunaga and Minami~\cite{matsunagaminami}\footnote{\ 
	though Hodgkin had already dubbed the map 
	``forgetful''~\cite[p.~72]{hodgkin1975kunneth}}
is stronger than 
the condition that $\KG(M;\Q) \lt \K(M;\Q)$
be surjective,
which Fok~\cite{fok2017formality} 
named \emph{rational K-theoretic equivariant formality}
and showed is equivalent to cohomological 
\emph{equivariant formality}
in the traditional sense~\cite{GKM1998} that the restriction 
$\HG(M;\Q) \lt \H(M;\Q)$ 
along the fiber inclusion in the Borel fibration
$M \to M_G \to BG$
be surjective.
Goertsches and Mare~\cite[Cor.~1.3]{goertschesmare2014}
showed a cohomogeneity-one action 
of a compact, connected Lie group $G$
on a smooth closed manifold $M$
with orbit space an interval 
is equivariantly formal if and only if $\rk G = \max\{\rk K^-,\rk K^+\}$,
so the same holds of rational K-theoretic equivariant formality
and the rank equation
is a necessary condition for K-theoretic equivariant
formality over the integers. 
The converse also holds, 
at least with the standard restriction on fundamental groups.

\eqf
\bpf
We consider the Hodgkin--K{\"u}nneth spectral sequence~\cite[Intro., Cor.~1, p.~6]{hodgkin1975kunneth} 
for the left multiplication $G$-action on $X = G$ and 
the given action on $Y = M$, 
a ($\Z \x \Z/2$)-graded left--half-plane spectral sequence
which starts at 
\[
E_2^{*,*} = \Tor^{*,*}_{RG}(K^*_G X, K^*_G Y) = 
\Tor^{*,*}_{RG}(\Z, K^*_G M)
\]
and,
given the hypothesis on $\pi_1 G$,
converges 
to 
\[K_G^*(X \x Y) = K_G^*(G \x M) \iso K^* (M)\mathrlap.\]
The forgetful map $K^*_G (M) \lt K^* (M)$ we wish to show is surjective
can be identified \cite[Prop.~9.1, p.~71]{hodgkin1975kunneth}
with the edge map 
\[
\KG (M) \longepi
\xt {RG}\Z{\KG (M)}
= E_2^{0,*} \longinc 
\Ei\col
\mathrlap.
\]
In each case we will
verify the groups
$\Tor^{\leq -1}_{RG}(\Z,K^*_G M) = 0$ vanish,
showing the spectral sequence collapses and the edge map is a surjection. 
We will repeatedly use the following facts.
First, if $K/H$ is an odd-dimensional sphere,
then $\rk K = 1 + \rk H$,
while if $K/H$ is an even-dimensional sphere,
then $\rk K = \rk H$.
Second~\cite[Thm.~3.6]{atiyahhirzebruch},
that for $\G$ closed and connected of full rank in $G$
we have $K^1(G/\G) = 0$
and $K^0(G/\G)$ free abelian (of rank $|WG|/|W\G|$).
Third~\cite[(7), p.~19]{gonzalezzibrowius2017SCSQ}, 
the groups $\Tor^{\leq p}_{RG}(\Z,R\G)$ vanish
for $\G \leq G$ closed and connected with $\rk G - \rk \G < |p|$,
so that particularly $\Tor^{\leq -2}_{RG}(\Z,R\G)$ 
vanishes for $\G \in \{\Kpm,H\}$.

\medskip

\nd\emph{Suppose $\rk G = \rk H + 1$}.

\smallskip
In these cases we know that one of $\Kpm$ has rank greater than that of $H$,
and our hypothesis on $\Kpm$ implies that $R\Kpm$ is polynomial~\cite{steinberg1975pittie},
so the corresponding restriction $R\Kpm \lt RH$ is surjective by 
\Cref{thm:Koddcircle,thm:Koddsphere} and 
the Mayer--Vietoris sequence of 
\Cref{thm:Kodd}
shows $K^1_G(M)$ vanishes,
leaving a short exact sequence of $RG$-modules
$K^0_G(M) \mono RK^- \x RK^+ \epi RH$. 
Applying the derived exact sequence of the functor $\Z \ox_{RG} {-}$,
we find $\Tor^{\leq -2}_{RG}(\Z,K^*_G M)$ vanishes as above.
Since in fact the $E_2$ page is only inhabited by 
$E_2^{0,0}$ and $E_2^{-1,0}$,
we know the former of these is $K^0(M)$ and the latter $K^1(M)$.
Thus the forgetful map will be surjective if and only if 
also $\Tor^{-1}_{RG}(\Z, K^0_G M) = K^1(M) = 0$.
Using the \MVS of the standard cover, 
we must show $K^0(G/K^-) \+ K^0(G/K^+) \lt K^0(G/H)$ is surjective
and $K^1(G/K^-) \+ K^1(G/K^+) \lt K^1(G/H)$ injective.

For surjectivity, 
assume without loss of generality that $\rk G = \rk K^+$,
so that $K^1(G/K^+)$ is zero
and $K^0(G/K^+)$ is free abelian;
in particular, then the \AHSS 
$\H(G/K^+) \implies \K(G/K^+)$ collapses.
There is an evident bundle map
\[
\xymatrix@C=.75em{
	K^+/H \ar[d]\ar[r]& {*}\ar[d]\\ 
	G/H  \ar[d]\ar[r]&  G/K^+\ar@{=}[d] \\ 
	G/K^+ \ar@{=}[r]&  G/K^+ 
}
\] 
inducing a map of {Atiyah--Hirzebruch--Leray--Serre spectral sequences}.
We have just seen the right spectral sequence collapses,
and the map then shows
all differentials out of the zero row of the left spectral sequence
must vanish as well. 
Particularly this means that the row $\Ei^{*,0}$ is a quotient of 
$E_2^{*, 0} = K^*(G/K^+)$; 
and since $K^+/H$ is an odd-dimensional sphere, 
$\K(K^+/H)$ is an exterior algebra $\ext[z]$ on one generator $z \in K^1(K^+/H)$,
so that
\[
E_2 = H^*\big(G/H;K^*(K^+/H)\big) \iso \H(G/H) \ox \ext [z].
\]
Since each diagonal thus contains
only one nonzero entry, we have $E_\infty \iso \K(G/H)$
as groups and thus, since odd columns are zero, $E_\infty^{*,0} \iso K^0(G/H)$.
This is a quotient of the row $E_2^{*, 0} \iso \H(G/K^+)$,
so the collapse $\H(G/K^+) \iso K^0(G/K^+)$ of the \AHSS on
the right shows $K^0(G/K^+) \lt K^0(G/H)$ is surjective.

Injectivity is obvious if $K^1(G/\Kpm) = 0$, 
so now assume as well that
$\rk K^- = \rk H = \rk G -1$.
We consider the map of Hodgkin--K{\"u}nneth spectral sequences 
corresponding to $X = G$ and $G/H = Y \to Y' = G/K^-$.
These are concentrated in the $0$-row
and again by the vanishing of $\Tor^{\leq -2}$, the spectral sequences 
both collapse at $E_2$, so the map $K^1(G/K^-) \lt K^1(G/H)$ 
may be identified with the map $\Tor^{-1}_{RG}(\Z,RK^-) \lt \Tor^{-1}_{RG}(\Z,RH)$.
But as $K^-/H$ is an even-dimensional sphere by assumption, 
\Cref{thm:Kevensphere} shows
$RH$ is free of rank two over $RK^-$,
so one has a short exact sequence $RK^- \mono RH \epi RK^-$.
Applying the derived exact sequence of $\Z \ox_{RG} -$
and the vanishing of $\Tor^{-2}$, we see
$\Tor^{-1}_{RG}(\Z,RK^-) \lt \Tor^{-1}_{RG}(\Z,RH)$
is injective as claimed. 

\medskip

\nd\emph{Suppose $\rk G = \rk H $}.

\smallskip

Since $K^1(G/\Kpm) = 0 = K^1(G/H)$ in this situation,
the sequence of \Cref{thm:mainK}
separates into the two short exact sequences
\[
0 \to K^0_G (M) \lt RK^- \x RK^+ \lt \defm B \to 0\mathrlap,
\]
\vspace{-1em}
\[
0 \to 	B \lt RH \lt K^1_G(M) \to 0
\]
of $RG$-modules.
From the vanishing of $\Tor^{\leq -2}$ we 
get $RG$-module isomorphisms 
\[
\phantom{\qquad(n \geq 1)\mathrlap,}
\Tor^{-n-2}_{RG}(\Z,K^1_G M) 
\iso 
\Tor^{-n-1}_{RG}(\Z,B) 
\iso 
\Tor^{-n}_{RG}(\Z,K^0_G M) 
\qquad(n \geq 1)\mathrlap,
\]
and from \Cref{thm:Keven} we also have an $RG$-module isomorphism
$K^0_G(M) \iso K^1_G(M)$, so the higher Tors are $2$-periodic.
But $\Z$ has finite projective dimension over $RG$ 
(indeed, the Koszul algebra $RG \ox \K G$ 
is a resolution of length $\rk G$),
so these higher Tors vanish. 
\epf

\brmk
The last sentence in this proof, the observation it concludes the proof, and
the request for such a result in the first place are all due to Marcus Zibrowius.
\ermk


{\footnotesize\bibliography{bibshort} }

\smallskip

\nd\footnotesize{%
\textsc{Department of Mathematics, 
		Imperial College London,
180 Queen's Gate,
London SW7 2AZ, UK\\
}
	\url{j.carlson@imperial.ac.uk}
}
\end{document}

%% file: tori.tex
\section{Coverings and mapping tori}\label{sec:mappingtori}

We begin with this section because
it is the only one involving any inversion 
of coefficients or any specifically equivariant homotopy theory.
It does not involve representation theory
or Lie theory in any serious way, so it is somewhat independent
of the rest of the document,
and we take it as an opportunity to get some long definitions out of the way.

Recall from
\Cref{thm:Mostert} 
that if a compact Lie group acts smoothly on a compact manifold $M$
with orbit space a circle (the case in \Cref{fig:torus}),
then $M$ is diffeomorphic to
the \defd{mapping torus} of the right translation 
by some element $w \in N_G(H)$ on $G/H$,
namely
\[
	\frac{G/H \x [0,1]}
			{(gH,1) \sim (gwH,0)}
		\mathrlap.
\]
As $w$ is of finite order $\defm{|w|}$, 
cutting the mapping torus at $t = 1$, 
gluing $|w|$ copies end to end, 
and then regluing the fiber $t = 0$ to $t = |w|$ by $w^{|w|} = \id_{G/H}$,
we see $G/H \x S^1$ is a $|w|$-sheeted covering of $M$.
The $G$-equivariant K-theory of $G/H \x S^1$ is easy to compute,
so most of our work is in computing the equivariant cohomology 
of a space from that of a finite-sheeted cover.

\bdefn\label{def:cohomology}
Let $G$ be a topological group.
A $\defm{G}$\defd{-}$\defm{n}$\defd{-cell} 
is a space $G/K \x D^n$,
where $K \leq G$ is a closed subgroup and $\defm{D^n}$ the closed $n$-disc,
equipped with the $G$-action $g\.(hK,x) \ceq (ghK,x)$.
A $\defm{G}$\defd{--CW complex} 
is a $G$-space $X$ constructed iteratively as the colimit (= union)
of a sequence of spaces $X_n$,
where $X_0$ is a disjoint union of $G$-$0$-cells
and otherwise each $X_n$ is obtained from $X_{n-1}$
by adjoining a collection of $G$-$n$-cells 
$G/K_\a \x D^{n}_\a$ along 
$G$-equivariant attaching maps $G/K_\a \x S^{n-1}_\a \lt X_{n-1}$.
When we do not specify otherwise, 
$\defm{S^n}$ comes equipped with the \emph{trivial} $G$-action
(and hence is, if you like, a $G$--CW complex
each of whose $G$-cells is of the form $G/G \x D^k$).
A $\defm G$\defd{--CW pair} $(X,A)$ comprises a $G$--CW complex $X$
and a $G$--CW subcomplex $A$, 
meaning each $G$-cell of $A$ is also a $G$-cell of $X$.
Given a $G$--space $X$, we denote by $\defm{X_{+}} \ceq X \amalg {*}$
the disjoint union of $X$ and a new isolated, $G$-fixed point $*$.

A \defd{reduced $G$-equivariant} ($\Z$-graded)
\defd{cohomology theory}
is a contravariant graded abelian group--valued homotopy functor
$\defm \tE = \Direct_{n \in \Z} \wt E^n$ 
on the category of pointed $G$--CW complexes
which takes a cofiber sequence
$A \to X \to X/A$
to an exact sequence of groups
and is equipped with a natural graded group isomorphism
$\defm\s\: \wt E^* X \isoto \wt E^{*+1} \susp X$
of degree one, the \defd{suspension},
where $\susp X = S^1 \wedge X$ is the reduced suspension of $X$. 
(Possibly obscured in the notation: 
$S^1$ is again assumed to have trivial $G$-action.)
Such a theory comes automatically with an associated 
\emph{unreduced theory} on unpointed $G$--CW pairs
given by $\defm{\E(X,A)} \ceq \tE(X/A)$ 
(by convention $X/\0 \ceq X_+$)
and satisfying the Eilenberg--Steenrod
axioms save dimension~\cite[\SS1]{matumoto1973equivariant}.

Let $\defm{\mathsf{Orb}_G}$ 
denote the category of orbits $G/K$ (for $K$ closed) 
and $G$-equivariant maps,
$\defm{h\Og}$ the category with the same objects 
but morphisms $G$-homotopy classes of $G$-maps,
$\defm\Top$ the category of topological spaces, and
$\defm\Ab$ the category of abelian groups.
A \defd{coefficient system} is a contravariant functor
$M\: h\Og \lt \Ab$.
For a given space $X$, 
the fixed point set assignment $G/H \lmt X^H$ 
gives a standard contravariant functor $\Og \lt \Top$
and composing any covariant functor $\Top \lt \Ab$
gives a coefficient system.
As an example, for each $n \in \N$
and each $G$--CW complex $X$
there is a functor
$
\defm{\ul H_n(X)}\: G/H \lmt H_n(X_n^H,X_{n-1}^H)
$.
The assignment $X \lmt \ul H_n(X)$ is itself covariantly functorial
in $G$--CW complexes.

The \defd{Bredon cohomology} 
$\defm{\HG(X;M)}$ of a $G$--CW complex $X$ 
with coefficients in a coefficient system $M$
is defined as the cohomology
of the complex 
$\defm{C^n_G(X;M)} \ceq \mr{Nat}\big(\smash{\ul H_n(X)},M\big)$
of natural transformations $\smash{\ul H_n} \lt M$,
where the $n\th$ coboundary map of the complex is 
precomposition with the tuple $\6_n = (\smash{\6_n^{G/H}})_{G/H \in \Og}$ 
for $\smash{\6_n^{G/H}}$ the connecting map
in the long exact homology sequence of the triple $(X_{n+1}^H,X_n^H,X_{n-1}^H)$.
Bredon cohomology is the unique unreduced $G$-equivariant cohomology theory
$\E$ which satisfies the wedge axiom 
and the requirement that $\E(G/H) = E^0(G/H) = M(G/H)$ for $G/H \in \Og$.

We write $\defm{|\G|}$ for the order of a group $\G$.
\edefn

\covering


\bpf
We first show the result for Bredon cohomology $H^p(-;E^q)$.
As the group $E^q(G/K)$ admits division by $|\G|$,
a classical 
Leray spectral sequence argument (apparently due to Grothendieck
\cite[Thm.~5.3.1, Cor. to Prop. 5.2.3]{grothendieck1957tohoku})
shows 
\[
\defm{\phi_{G/K}}\:\H\big(X^K_p/\G,X_{p-1}^K/\G; E^q(G/K)\big)
\lt \H\big(X_p^K,X_{p-1}^K; E^q(G/K)\big){}^\G
\]
is an isomorphism. 
Endow $E^q(G/K)$ with the trivial $\G$-action.
Since the Kronecker pairing is $\G$-invariant,
the universal coefficient morphism 
\[
H^p\big(X_p^K,X_{p-1}^K;E^q(G/K)\big) \longepi \Hom\big(H_p(X_p^K,X_{p-1}^K),E^q(G/K)\big)
\]
is also $\G$-equivariant,
and since $E^q(G/K)$ is divisible by $|\G|$,
induces a surjection of $\G$-invariants, 
every $\G$-invariant element being the average over a $\G$-orbit.
It follows from this surjectivity, 
the surjectivity of $\phi_{G/K}$,
and the functoriality of
the universal coefficient theorem that
\[
\defm{{f_{G/K}}}\:\smash{\Hom\!\big(\mnn H_p(X_p^K/\G,X_{p-1}^K/\G),E^q(G/K)\big)
	\lt \Hom\!\big(H_p(X_p^K,X_{p-1}^K),E^q(G/K)\big)^\G}
\]
is also a surjection.
By the observation that 
$\big(\frac{G \x \G}{H \x \D}\big){}^K/\G
		\homeo
		\Big(\big(\frac{G \x \G}{H \x \D}\big)/\G\Big){\vphantom{\big)\mn}}^K$\!,
our assumption on the isotropy groups of $X$,
and induction,
we have $(X_n)^K/\G = (X_n/\G)^K$ for all $n$,
so the natural transformations $\ul H_p(X/\G) \lt E^q$
are encoded by coherent sequences in the domain 
of $\prod_{G/K \in \Og}  f_{G/K}$. 
Equally, assigning each $E^q(G/K)$ the trivial $\G$-action,
the $\G$-equivariant natural transformations
$\ul H_p(X) \lt E^q$
are coherent sequences in the codomain of 
$\prod_{G/K \in \Og}  f_{G/K}$.
Thus we will have an isomorphism
$C^p_G\big(X/\G;E^q) \isoto C^p_G(X;E^q)^\G$
if we can show $ f_{G/K}$ is also injective for each $G/K \in \Og$.

To this end we may forget the corestriction to $\G$-invariants
in the codomain
and just show the map of Homs is injective,
and for this it is enough to see the predual
\[
	\defm{\psi_{G/K}}\:
H_p(X_p^K,X_{p-1}^K) 
	\lt
H_p\big((X_p/\G)^K,(X_{p-1}/\G)^K\big) 
\]
is surjective.
From the definition of a $(G \x \G)$--CW complex
and our assumption on isotropy groups,
the quotient $X_p^K/X^K_{p-1} = (X_p/X_{p-1})^K$
is a wedge of summands
\[
(G/H_\a \x \G/\D_\a)^K_+ \^ S^p = 
\big((G/H_\a)^K \x \G/\D_\a\big){}_+ \^ S^p
\] 
for various product subgroups $H_\a \x \D_\a \leq G \x \G$,
so the group 
$H_p(X_p^K,X^K_{p-1}) \iso \wt H_p (X_p^K/X^K_{p-1})$
decomposes as
\[
\Direct_\a \wt H_p\Big(\mn\big((G/H_\a)^K \x \G/\D_\a\big)_+ \^ S^p\Big)
\iso
\Direct_\a \wt H_0\big((G/H_\a)^K \x \G/\D_\a\big)_+
\iso
\Direct_\a  \H_0\big(\mnn(G/H_\a)^K\big)^{\oplus\,|\G/{\D_\a}|}
	\mathrlap,
\]
and quotienting by $\G$ we have a similar isomorphism
\[
H_p\big((X_p/\G)^K,(X_{p-1}/\G)^K\big)
	\iso 
\smash{\Direct_\a H_0\big(\mnn(G/H_\a)^K\big)}
	\mathrlap.
\]
But under these identifications 
the $\a\th$ summand of $\psi_{G/K}$
is just iterated addition
$(x_1,\ldots,x_{|\G/\D_\a|}) \lmt x_1 + \cdots + x_{|\G/\D_\a|}$
in the group $H_0\big(\mnn(G/H_\a)^K\big)$,
which is certainly surjective.

Varying $p$, we have our isomorphism of cochain complexes
$C^*_G(X/\G;\E) \lt C^*_G(X;\E)^\G$.
Note that $C^*_G(X;\E)$ is divisible by $|\G|$ 
and recall that given a cochain complex $C$ 
of $|\G|$-divisible $\G$-modules,
the inclusion $C^\G \longinc C$
induces an isomorphism $\H(C^\G) \isoto \H(C)^\G$
and multiplication by $|\G|$ is again invertible on $\H(C)$.
Finally the composite
\[
H^*_G(X/\G;\E) \isoto \H\big(C^*(X;\E)^\G\big) \isoto H^*(X;\E)^\G
\]
is the claimed isomorphism in Bredon cohomology.

There is a equivariant \AHSS
due to Matumoto~\cite[\SS4]{matumoto1973equivariant},\footnote{\ 
	The spectral sequence with sheaf coefficients due to Segal%
	~\cite[\SS5]{segal1968equivariant}
	reduces to this one in the case $\E = \KG$ but is less immediately
	adapted to our needs.
}
functorial in and converging to the $\E$-cohomology of
finite $G$--CW complexes, 
and the entries $E_2^{p,q}$ of its second page are 
the Bredon cohomology groups
$H_G^p(-;E^q)$ with coefficients in the coefficient system 
$K \lmt E^q(G/K)$.
Forgetting the $\G$-action and regarding $X$ as a $G$--CW complex,
we see $\pi\: X \lt X/\G$ induces a morphism of these spectral sequences.
Since the spectral sequence can be defined using a Cartan--Eilenberg
$H(p,q)$-system with $H(p,q) \ceq \Direct_n E^n(X_{p-1},X_{q-1})$ 
and the skeleta $X_j$ are $\G$-invariant by definition,
the differentials $d_r$ of this spectral sequence are $\G$-equivariant.
On $E_2$ pages, the induced map of spectral sequences is 
$\HG(X/G;\E) \lt \HG(X;\E)$,
which we have just seen is an isomorphism onto its image $\HG(X;\E)^\G$.
%
Inductively applying the recollection 
about invariants of cochain complexes from the previous paragraph to each page,
we see $\pi^*$ induces a pagewise isomorphism 
of one spectral sequence with the $\G$-invariants of the second,
and so at $E_\infty$ we recover an isomorphism
$\gr \E(X/\G) \simto (\gr \E X)^\G$,
where $\defm\gr$ denotes the associated graded module with respect
to the cellular filtration.
But for any filtered $\G$-module $N$ divisible by $|\G|$,
the inclusion $N^\G \longinc N$
induces an isomorphism $\gr (N^\G) \isoto (\gr N)^\G$,
so the $E_\infty$ map further factors through an isomorphism
$\gr \E(X/\G) \isoto \gr\big(\E (X)^\G\big)$.
This is the associated graded map induced by $\E(X/\G) \lt \E(X)^\G$,
so as the filtration involved is finite,
that map is an isomorphism as well%
~\cite[Thm.~2.6]{boardman1999conditionally}.
\epf

As a corollary we have a result on mapping tori,
which we prefer to state as a ring isomorphism,
so we will need to define an additional notion.

\bdefn\label{def:multiplicative}
A $G$-equivariant cohomology theory $\E$ is said to be 
\defd{multiplicative} if
$\E$ is valued in commutative graded algebras 
and the suspension axiom is replaced in the following way.
Note that $\E(*,\0) = \wt E^0 S^0$ 
is a commutative ring with unity $1$
and the projections
$\pi_Y,\pi_X\: Y \x X \lt Y,X$
induce a natural \emph{cross product}
\eqn{
	\tE Y \ox \tE X 	&\os{\defm\x}\lt \tE(Y \^ X),\\
	y \ox x		&\lmt \pi_Y^* y \. \pi_X^* x.
}
The new axiom is that there exist an element 
$\defm{\es} \in \wt E^1 S^1$
such that the map
$\s\: \tE X \isoto \wt E^{*+1} (S^1 \^ X)$ 
given by $\s(x) \ceq \es \x x$
is a natural isomorphism.
\edefn

\brmk
This is somewhat leaner than the usual axiomatization.
It is typical in defining a multiplicative cohomology theory
to demand it be represented by a ring spectum,
but we do not require our theories to satisfy the wedge axiom,
and thus our results will allow 
for things like $p$-completed theories.

For non-represented theories, it is usual 
to require natural cross products satisfying naturality 
axioms, but it seems simpler to demand cup products
and instead note the other axioms follow from the \CGA structure
and functoriality.
The typical axiomization also demands 
sign-commutativity of evident squares involving suspensions, 
but these are all consequences of 
graded commutativity 
and the uniform definition of suspension
as a cross product.
Unreduced theories additionally require the cross product 
cooperate with the connecting maps from 
the long exact sequences of a pair, 
but the connecting map can be defined in terms of the suspension 
in the unreduced theory,
so the commutativity of these squares is again a formal consequence
of functoriality and the uniform definition of the suspension.
\ermk

Now we can state the result.

\begin{restatable}{lemma}{mappingtorus}\label{thm:mappingtorus}
	Let $Y$ be a $G$-space and $\varphi$ a self-homeomorphism 
	of $Y$ commuting with the $G$-action
	and such that there exists a positive integer $\defm \ell$
	such that $\varphi^{\ell}$ is homotopic to $\id_Y$.
	Write $\spac$ for the mapping torus of $\varphi$ and
	let $\E$ be a $\Z$-graded multiplicative equivariant cohomology theory
	valued in $\Z[\sfrac 1{\,\ell}]$-algebras.
	Write $\defm{\E} \ceq \E\mn\,\mn(*)$.
	Then
	\[
	\E \spac \iso \xt{\E}{\E(Y)^{\ang{\varphi^*}}}{\ext_{\E}[z]}
	\mathrlap,
	\] 
	where $z$ is the pullback 
	of a generator of $\wt E^1(S^1) \iso \wt E^0(S^0) = \E$
	under $X \lt S^1$.
\end{restatable}
	
	Here, as usual, $\defm{\E(Y)^{\ang{\vp^*}}}$ 
	denotes the subring of elements
	invariant under pullback by $\vp$.

\bpf
Note that $\spac$ admits an $\ell$-sheeted cyclic covering $\defm Z$
by the mapping torus of $\varphi^\ell$,
which is homeomorphic to the mapping torus $Y \x S^1$
of the identity. 
This homeomorphism takes the covering action to a $\Z/\ell$-action on 
$Y \x S^1$ under which $1 + \ell\Z$ acts, up to homotopy,
as $(y,\t) \lmt \smash{\big(\varphi(y),\t + \frac {2\pi}\ell}\big)$,
which, rotating the $S^1$ component,
is in turn homotopic to $(y,\t) \lmt \smash{\big(\varphi(y),\t\big)}$.
It follows from the suspension axiom 
for $\wt E^*$ that 
$\E S^1 					\iso 
 \E \+ \smash{\wt E^* S^1} 	\iso 
 \E \+ \E[1] 				\iso 
 {\E} \+ \E\mn \. \mn\{z\}$.
Now assuming multiplicativity,
as $z \in E^1 S^1$ is a free $E^0({*})$-module generator of $\wt E^* S^1$,
we have $\E S^1 \iso \ext_{\E}[z]$.
It follows again from the suspension axiom that
$\E S^1 \ox_{\E}\E Y \lt \E(S^1 \x Y)$
is a ring isomorphism.\footnote{\ 
	Explicitly, naturality of multiplication implies
	the suspension isomorphism 
	$\wt E^* (Y_+) \lt \wt E^{*+1}(S^1 \wedge Y_+)$
	is given by multiplication by the pullback of $z$,
	giving a natural nonunital ring isomorphism
	$\wt E^* S^1 \ox_{\E} \wt E^* (Y_+) 
	\lt 
	\wt E^*(S^1 \wedge Y_+)$.
	From the cofiber sequence 
	$S^1 \vee Y_+ \to S^1 \x Y_+ \to S^1 \wedge Y_+$
	we get
	$\E S^1 \ox_{\E}\E (Y_+) \isoto \E(S^1 \x Y_+)$
	and from
	$Y \to Y_+ \leftrightarrows *$
	we get
	$\E S^1 \ox_{\E}\E Y \isoto \E(S^1 \x Y)$.
}
The action of $1+\ell\Z$
on $\E Y \ox_{\E} \E S^1 \iso \E Z$
is given by $a \ox s \lmt \varphi^* a \ox s$,
so an application of \Cref{thm:covering}
yields the claim.
%
\epf

\bprop\label{thm:circle}
Let a cohomogeneity-one action of a compact, connected 
Lie group $G$ on a smooth manifold $\man$ be 
given with orbit space $\man/G \homeo S^1$.
Recall from \Cref{thm:Mostert} that this means
$\man$ is $G$-equivariantly diffeomorphic 
to the mapping torus of 
right multiplication on $G/H$
by some element $w\in N_G(H)$
and let $\ell$ be the smallest 
positive integer such that $w^\ell$
lies in the identity component of $N_G(K)$.
Suppose 
$\E$ is a $\Z$-graded multiplicative equivariant cohomology theory
valued in $\Z[\sfrac 1{\,\ell}]$-algebras.
Then one has a graded ring isomorphism
\[
\E \man \,\iso\, \xt{\E}{\E(G/H)^{\ang{r_w^*}}} {\ext_{\E}[z_1]}, \quad |z_1| = 1.
\]

\eprop
\bpf

Note that $w^\ell$ lies in the path-component of the identity,
so that right multiplication by $w^\ell$ is homotopic to $\id_{G/H}$, 
and apply 
\Cref{thm:mappingtorus}.
%
\epf

The result we want follows immediately:

\Kcircle
\brmk
There is a transfer map in K-theory we could also apply directly
to bypass this level of generality.
\ermk

Such a clean statement is not possible without inverting the 
order of $w$.

\bex\label{ex:inversion}
Let $G = \SO(n)$ and $K$ the block-diagonal subgroup 
$[1]^{\+ n-2} \+ \SO(2)$. 
Then $N_G(K)$ has two components, 
represented by the identity matrix and the block-diagonal
$w = [1]^{\+ n-3} \+ [-1] \+	\bigl[\begin{smallmatrix}
0&1\\1&0 
\end{smallmatrix}\bigr]$,
conjugation by which corresponds to complex
conjugation under the standard identification
of $\U(1)$ with the unit circle in the complex plane.
Thus $w$ acts on $R\SO(2) \iso \Z[t^{\pm 1}]$ 
by $t \bij t\-$, where $t\: \SO(2) \isoto \U(1)$ 
is the defining representation on $\C \iso \R^2$.
We let $M$ be the mapping torus of the right action 
of ${w}$ on $G/K$.
To proceed integrally rather than via \Cref{thm:Kcircle},
we use the \MVS of the cover of $M$ by two intervals
overlapping at the endpoints.
This is an exact sequence
\[
0 \to K^0_G M \lt RK \x RK \lt RK \x RK \lt K^1_G M \to 0
\]
where the middle map is $(a,b) \lmt (a-b,a-wb)$.
Since the first map is diagonal, 
the middle map may be replaced with the map
$\phi\: RK \lt RK$ taking $a$ to $a - wa$.
Thus 
\[
K^0_G(M) \iso \ker \phi = R(K)^{\ang w} = \Z[t+t\-]\mathrlap,
\] 
\[
K^1_G(M) \iso \coker\phi \,=\, \quotientmed{\Z[t^{\pm 1}]\,}{\Z\{t^n - t^{-n} : n \in \N\}}	\mathrlap.
\]
Since the denominator in the cokernel induces on the numerator 
precisely the relations $t^{-n} \equiv t^n$,
a set of coset representatives for $\coker \phi$ is given by $\Z\{1,t,t^2,t^3,\ldots\}$.
Writing $q = t + t\-$, one sees
\[
[1] \os q\lmt [t]+[t\-] = 2[t],\qquad
[t] \os q\lmt [t^2 + 1] 
\os q\lmt [t^3 + 3t] 
\os q\lmt [t^4 + 4t^2 + 3] 
\os q\lmt \cdots,
\]
and generally $q^n \. [t]$ has highest term $[t^{n+1}]$,
so $K^1_G\big(M;\Z[\sfrac 1{\,2}]\,\mn\big)$ is a free cyclic 
$K^0_G\big(M;\Z[\sfrac 1{\,2}]\big)$-module on $[1]$. 
Note that with $\Z$ coefficients, 
$K^1_G(M)$ is not a free $K^0_G(M)$-module.
\eex

%% file: MVgeneral.tex

\MV

The additional structure on the connecting map
is most helpful when even or odd cohomology
of the constituent subsets vanishes,
making the connecting map surjective.

\bex\label{ex:3mf}
Let $M$ be a closed, oriented $3$-manifold.
Then $M$ can be triangulated.
A regular neighborhood $U$ of its $1$-skeleton
is an open handlebody (i.e., 
homeomorphic to the bounded component cut out
of $\R^3$ by an embedded closed surface),
and examining the local picture in 
each $3$-simplex, 
one sees the interior of the complement $V$ is
also a handlebody. 
The closures of $U$ and $V$ meet in 
a closed, oriented surface $S_g$,
and this assemblage is called a \emph{Heegaard splitting}
of $M$.
Letting $N_g$ denote a standard genus-$g$ handlebody
with boundary $S_g$,
we may write $M \homeo N_g \union_f N_g$
for some gluing homeomorphism $f\: S_g \lt S_g$.
If we write $a_j$ for the standard 
$g$ circles generating $H_1(N_g)$
and $b_j$ for the $g$ circles bounding discs in $S_g$
representing the other standard generators,
so that $|a_i \inter b_j| = \smash{\delta^{i}_{j}}$,
then $M$ is determined up to homeomorphism 
by the images $f(b_j)$.
Let $\a_j$ and $\b_j$ be the dual basis of $H^1(S_g)$.

Fattening $U$ and $V$ slightly, 
we may apply the \MVS in cohomology,
which contains the subsequence
\[
0 	\to H^1(M) 
	\os\vk\lt \Z^g \+ \Z^g
	\os\l\lt H^1(S_g)
	\os\d\lt H^2(M) 
	\to 0\mathrlap.
\]
Thus $H^1(M)$ and $H^2(M)$ are determined by the map $\l$,
which is in turn determined by the map $f$.
If we make the identifications 
$U \inter V = S_g \subn N_g = U$,
then the first component $\l_1\:\Z^g \lt H^1(S_g)$
is the inclusion $\iota^*\: \a_j \lmt \a_j$
and the second component $\l_2$ is $f^*\iota^*$,
so we have an isomorphism
\[\im \d \iso \coker \l = \frac{\Z\{\a_j,\b_j\}}{\Z\{\a_j,f^*\a_j\}}\mathrlap,\]
which in particular is spanned by the images
of the $\b_j$,
and $H^1(M) \iso \ker \l$
is spanned by elements $(\sum m_i \a_i,\sum n_j \a_j)$
such that $\sum n_j f^*\a_j$ has no $\b$-component.
By \Cref{thm:MV}, 
the cup product $\mu_{1,2}\: H^1(M) \x H^2(M) \lt H^3(M)$
is determined by 
$y \cup \d(z) = \d(\l_1 \vk y \cup z)$,
where 
$\l_1 \vk y$ is some linear combination of the $\a_i$
and $z$ can be taken to be a linear combination of the $\b_j$,
and the second cup product is taken in $\H(S_g)$.
Since this product
is given on generators by $\a_i \cup \b_j = \d^i_j$,
the \MVS gives $\mu_{1,2}$ in terms of $H_1(f)$.
\eex

Though \Cref{thm:MV} 
does not seem to appear as such in the literature,
with a bit of faith 
it is possible to cobble together a proof from citations.

\begin{proof}[Terse proof of \Cref{thm:MV}]
	In the long exact sequence of a pair $(X,A)$,
	the connecting map $\E(A) \lt E^{*+1}(X,A)$
	is an $\E(X)$-module homomorphism;
	see Whitehead~\cite[(6.19), p.~263]{whitehead1962generalized}
	for an algebraic proof for cohomology theories represented by ring 
	spectra and note the proof still follows from our axioms.
	Up to homotopy, the Mayer--Vietoris sequence of $(X;U,V)$
	is the long exact sequence of a pair $(X', U' \amalg V')$
	in which $X'$ is homotopy equivalent to $X$ via
	a homotopy equivalence $X' \lt X$
	sending disjoint $G$--CW subcomplexes $U'$ and $V'$ 
	respectively to $U$ and $V$;
	\emph{cf.} Adams~\cite[p.~213]{adamsstable} for a 
	version of this statement for a representable theory.\footnote{\ 
		 Another version of this statement 
		appears in a MathOverflow solution
		due to J. Peter May~\cite{MO:MV} for CW-spectra
		(or, to quote, ``any halfway reasonable category''
		of spectra).
	}
\end{proof}

This in a moral sense a geometry paper, 
so for those with less faith, 
a more expansive and geometric account follows.

\begin{notation}
In what follows between now and the return to K-theory, 
all maps will be equivariant
with respect to a fixed topological group $G$
and all $G$-spaces will
come equipped with a $G$-fixed basepoint $\defm{*}$.
The wedge sum and smash product inherit the expected actions,
and the closed unit interval $I = [0,1]$
and circle $S^1 = I / (0 \sim\mspace{-1mu}1)$
are basepointed at $0$
and equipped with the trivial $G$-action.
We write $\defm{CX} = I \^ X $ for the reduced cone
and $\defm{\susp X} = CX/X =  S^1 \^ X$ for the reduced suspension,
with the induced actions.
\end{notation}


%

The $G$-structure is just along for the ride
in the proof that follows,
and everything we state through to 
\Cref{thm:oddzero} follows
for nonequivariant theories 
through the expedient of setting $G = 1$.

\defn
Let $\tE$ be a multiplicative $G$-equivariant cohomology theory
(not even necessarily equipped with suspension maps).
The diagonal
$	\defm{\Delta}\: 	X \lt X \^ X$
makes a $G$-space $X$ a \emph{coalgebra} in the sense that 
$
(\Delta \^ \id)\o \Delta =  
(\id \^\, \Delta)\o \Delta
$.
A right $X$-\emph{coaction} 
$\defm{\Delta_Y}\: Y \lt Y \^ X$ on a $G$-space $Y$ 
is a map such that
$
(\D_Y \^\, \id)\o \Delta_Y = (\id \^\, \D) \o \Delta_Y
$;
such a map makes $Y$ a right $X$-\emph{comodule}
and induces an additive homomorphism
$\Delta_Y\o \mu_{Y,X}\: \tE Y \ox \tE X \lt \tE Y$
which one checks, unravelling definitions,
to be a right $\tE X$-algebra structure. 
A map $f\: Y \lt Z$ between right $X$-comodules such that
$\Delta_Z \o f = (f \^ \id) \o \Delta_Y$ is an \emph{$X$-comodule homomorphism},
and induces a $\tE X$-algebra homomorphism $f^*\: \tE Z \lt \tE Y$.

\bprop\label{thm:pair}
Let $G$ be a topological group and
$\E$ a multiplicative $G$-equivariant cohomology theory.
Then in the long exact sequence of a $G$--CW pair $(X,A)$,
all objects are $\E X$-modules and all arrows
$\E X$-module homomorphisms.
In particular the image of $\E(X/A) \lt \E X$ is an ideal
and the image of $\E A \lt \wt E^{*+1} (X/A)$ is a nonunital
subring with zero multiplication.
\eprop

We adapt a proof 
from Hatcher's manuscript K-theory text~\cite[Prop.~2.15]{VBKT},
which considers the cross product with a single element
and does not make explicit use of the notion of a comodule.

\bpf
It will be enough to prove the result for the reduced theory $\wt E^*$.
Note that for pointed $G$--CW subcomplexes $A$ of $X$
and pointed $G$--CW complexes $S$ with trivial action,
$S \^ A$ admits the $X$-coaction
$s \^ a \lmt s \^ a \^ a$
and $S \^ (X \union CA)$ the $X$-coaction
\eqn{
	\phantom{a \^}	s 	\^ x	&\lmt  
	s \^ x \^ x,\\
	s \^ t \^ a	
	&\lmt   s \^ t \^ a \^ a.
}
It is easy to check these coactions make a cofiber sequence
$A \to X \to X \union CA$
a sequence of $X$-comodule homomorphisms.
To see this also makes the Puppe sequence
\[
A \os i\lt X \lt X \union CA \lt \susp A \os{\susp i}\lt \susp X \lt 
\susp (X \union CA) \lt \susp^2 A \lt \cdots
\]
a sequence of $X$-comodule homomorphisms,
it suffices to observe the coaction commutes with (suspensions of) the 
connecting map
$X \union CA \lt S^1 \^ A$
given by $t \^ a \lmt (1-t) \^ a$ and $x \lmt *$.
To replace $S\^ (X \union CA)$ with $S \^ X/A$,
observe the coaction $s\^[x]\lmt  s \^ [x] \^ x$
on the latter makes the collapse map another $X$-comodule homomorphism.

Applying $\tE$ to the Puppe sequence then yields an
$\tE X$-module structure on the long exact sequence of $(X,A)$.
To see the image of the connecting map has trivial multiplication,
note this map can be written as $\tE \susp A \lt \tE (X/A)$.
\epf

\brmk
The meticulous reader will observe 
that the proof of \Cref{thm:pair}
makes use of the fact
the coaction smashes with $X$ on one side and
the suspension smashes with $S^1$ on the other.
This choice actually matters;
the choice of a left $\E X$-action instead of a right 
requires an additional sign, making the connecting map
fail to be an $\E X$-module homomorphism.\footnote{\ 
In detail, for singular cohomology, 
the $\k$-submodule $C^*(X,A;\k)$ of cochains vanishing 
on $C_*(A)$ is a two-sided ideal of $C^*(X;\k)$ 
with respect to the cup product,
which thus restricts to both a right and a left action
of $C^*(X;\k)$ on $C^*(X,A;\k)$.
Using the zig-zag lemma to compute the connecting map $\d$ of the short exact sequence $C^*(X,A;\k) \to C^*(X;\k) \to C^*(A;\k)$ 
of cochain complexes gives 
$\d\big(a \cup i^*(x)\big) = \d a \. x$ but
$\d\big(i^*(x)\cup a\big) = (-1)^{|x|}x\.\d a$.
In terms of our preceding discussion, 
the sign arises because 
the connecting map of the pair $(X,A)$
factors as the composition of ring homomorphisms and the 
suspension isomorphism 
$\smash{H^*(A;\k) \os\d\to H^{*+1}(CA,A) \os\sim\from \wt H^{*+1}(\susp A)}$
arising from the long exact sequence of the pair $(CA,A)$
and the standard homeomorphism $CA/A \homeo \susp A$;
but since the suspension isomorphism  
can be identified as
$\smash{H^*(A;\k) \simto H^1 (S^1;\k) \ox_k H^*(A;\k) 
	\os\x\to \wt H^{*+1}(S^1 \^ A)}$,
the cross product \emph{on the left} with the fundamental class of $S^1$,
a sign can be avoided only by switching the side on 
which $H^*(X;\k)$ acts.
}
One could be forgiven for suspecting this has something to do with 
the well-known 
sign in the Puppe sequence:
our choice of $q\: t \^ a \lmt (1-t) \^ a$ 
for the map $X \union CA \lt A \^ S^1$
comes from a nonstandard identification 
$\smash{CX \union CA \epi \susp A \os -\to \susp A}$
in transitioning from the iterated cofiber sequence
to the Puppe sequence.
This choice of identification makes $\E q$ the \emph{opposite} $-\d$
of the connecting map $\d\: \E A \lt E^{*+1}(X,A)$
defined through the axioms
but makes the next map $\susp E^* i$ rather than the $- \susp E^* i$
it would become under the standard identification.
As $q$ and its variant $-q$ are both $X$-comodule maps, 
the choice between them is immaterial to the success of \Cref{thm:pair},
and moreover,
this choice inflicts a global sign of $-1$ on the connecting maps 
in each degree, so the correction factor 
arising from putting the $\E X$-action on the left 
would be a separate, logically independent sign. 
\ermk

To obtain the same result on connecting maps 
for the Mayer--Vietoris sequence,
we realize it as the long exact sequence of a pair,
as in the terse proof.

\begin{figure}[H]
	\caption{Schematic of $CU \union X' \union CV$ in \Cref{thm:pairMV}}%
	\label{fig:MV}%
	\centering{
		\includegraphics[height=4.2cm]{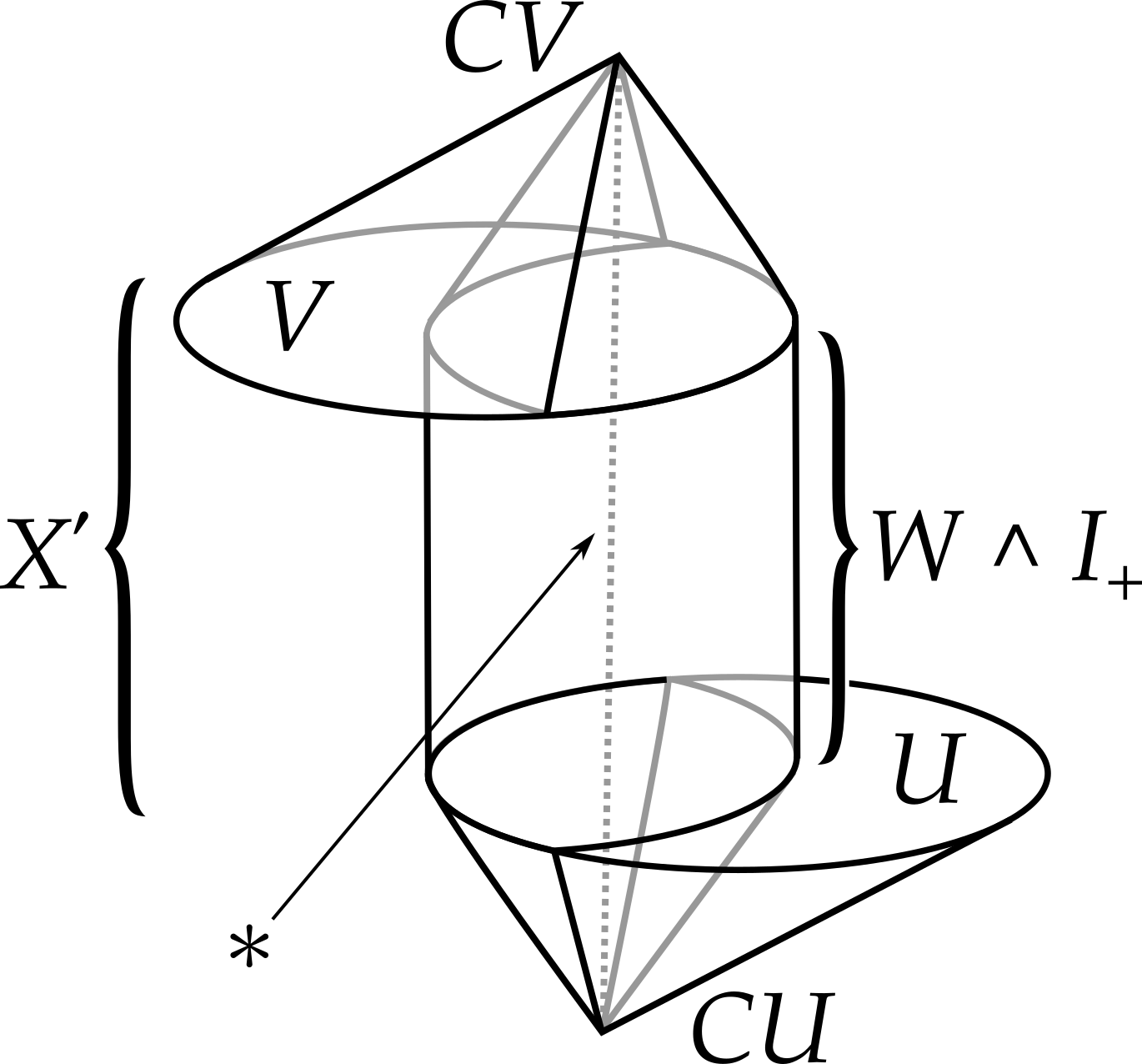}
	}
\end{figure}

\bprop\label{thm:pairMV}
Let $(X;U,V)$ be a triad of $G$--CW complexes with $X = U \union V$.
Write $W$ for the intersection $U \cap V$
and $\defm{X'}$ for the double mapping cylinder
$\big(U \x \{0\}\big) 
\union (W \x I) 
\union 
\big(V \x \{1\}\big)$ 
of the inclusions $U \hookleftarrow W \inc V$.
Then for any $G$-equivariant cohomology theory,
the long exact sequence of the pair
$\big(X', U \x\{0\} \amalg \mspace{-1mu}V \x \{1\}\big)$
is the Mayer--Vietoris sequence of the triad $(X;U,V)$.
\eprop
%
%

\bpf
It is again enough to assume $W$ is pointed 
and prove the result for the reduced theory.
In so doing, we replace $W \x I$ with the reduced cylinder
$W \^ I_+ = (W \x I) / \big({\{*\}} \x I\big)$,
turning $X'$ into $X'' = X' / \big({\{*\}} \x I\big)$
and $U \x \{0\}\amalg \mspace{-1mu}V\x\{1\}$ into $U \vee V$, 
which is naturally a subspace of $X''$
since the basepoints $(*,0)$ and $(*,1)$ have been identified.
The result is as in \Cref{fig:MV}.

Note $X''$
is $G$--homotopy equivalent to $X$ via the map collapsing the 
$I$-direction in the reduced cylinder $W \^ I_+$.
The Puppe sequence begins
\[
U \vee V \lt X'' \longinc CU \union X'' \union CV \xrightarrow{\, /X'' } \susp U \vee \susp V \lt \susp X''.
\]
We can replace the third term with $\susp W$
because the map collapsing $CU \vee CV$ to a point is a $G$--homotopy equivalence.
The maps then yield an exact sequence of graded groups
\[
\tE U \+ \tE V 
\longfrom
\tE X 
\os{\defm \d}\longfrom 
\wt E^{*-1} W 
\os{\defm \z}\longfrom
\wt E^{*-1} U \+ \wt E^{*-1} V
\longfrom \wt E^{*-1}X,
\]
which we check is the Mayer--Vietoris sequence:
\bitem
\item
That $U \vee V \longinc X''$ yields the pair of restrictions 
$\tE X \lt \tE U \+ \tE V$
is clear.
\item
The connecting map in the \MVS
is defined as the composition
\[
\wt E^{*-1} W \lt 
\tE(V/W) \xleftarrow{\ \sim\ } 
\tE(X/U) \lt 
\tE X,
\]
where the first map is the connecting map in the long exact sequence
of the pair $(V,W)$,
hence induced by $V/W \os\sim\from V \union CW \to \susp W$, 
the second is the excision arising 
from the homeomorphism $V/W \lt X/U$,
and the last is induced by the projection 
$X \longepi X/U$.
Thus the Mayer--Vietoris connecting map is obtained by following 
the path from $X$ to $\Sigma W$ along the 
bottom of the following commutative diagram,
while $\d$ comes from following along the top:
\[
\xymatrix@R=4.5em{
	\susp W & 
	X'' \ar@{->>}[l]_(.43){/ U \vee V} \ar@{->>}[r]^\sim \ar@{->>}[dl]^{/U}& 
	X 
	\ar@{->>}[d]
	\\
	V \union CW 
	\ar@{->>}[u]^{/V}
	\ar@{->>}[r]^(.55)\sim_(.53){/CW}
	& V/W \ar[r]^{\homeo} & X/U	.
}
\]
\item
The map $\z$ is induced as the composition along the right
in the commutative diagram
\[
\xymatrix@R=4em@C=1em{
	\susp W \ar[r]^(.24)\homeo&
	CW \union (W \wedge I_+) \union CW  
	\ar@{->>}[d]^{/ (W \^ I_+)} 
	\ar@{}[r]^(.47){}="a"^(.57){}="b"
	\ar@{^{(}->}"a";"b"&  
	CU \union (W \wedge I_+) \union CV
	\ar@{->>}[d]^{/ (W \^ I_+)} \\
	&
	\susp W \vee \susp W 
	\ar@{}[r]^(.21){}="a"^(.8){}="b"
	\ar@{^{(}->}"a";"b"& 
	\susp U \vee \susp V.
}
\]
On the other hand,
the left vertical map
collapsing a cylinder's worth of $W$s
is $G$--homotopy equivalent to the pinch map 
$\susp W \longepi \susp W \vee \susp W$
collapsing 
only the equator $W \x \{1/2\}$,
so the composition 
$\susp W \to \susp W \vee \susp W \to \susp U \vee \susp V$
is homotopic to $\susp j_U - \susp j_V$,
where $j_U,j_V\: W \longinc U,V$ are the inclusions.
The minus sign comes from observing 
a small neighborhood the cone point of the abstract $CU = U \wedge I$ 
lies near suspension coordinate $t=0$, 
agreeing with the suspension coordinate of the included copy of $CU$ in 
$CU \union (W \wedge I_+) \union CV$,
while the cone point of the included copy of $CV$ is near $t = 1$,
disagreeing with that of the abstract $CV$.
\qedhere
\eitem
\epf

The conjunction of these two results gives \Cref{thm:MV}.
Taking $W = U \inter V$ in the statement,
the image of $\d\: E^{*-1} W \lt \E X$
is an ideal with multiplication zero,
since $\d$ is induced by 
$X \lt \susp W$ and the multiplication of
the non-unital algebra 
$\tE\susp W$ is zero.
This result allows us to completely compute the ring $\E X$
from $\E U$, $\E V$, and $\E W$ in amenable cases.
We write $j_U,j_V\: W \longinc U,V$ and $i_U,i_V\: U,V \longinc X$.

\bprop\label{thm:oddzero}
Let $\E$ be a $\Z$-graded $G$-equivariant multiplicative cohomology theory
and $(X;U,V)$ a triple of $G$--CW complexes with $X = U \union V$
and such that the odd-dimensional $E$-cohomology of $U$, $V$, and $W = U \cap V$
vanishes. 
Then one has a graded ring and a graded $\E W$-module isomorphism,
respectively:
\[
E^{\mr{even}} X \,\iso\, \xu{\E W}{\E U}{\E V},\qquad
E^{\mr{odd}}  X \,\iso\, 
\Big(\mn\quotientmed{\E W\,}{\,\im j_U^* + \im j_V^*}\Big)[1].
\]
The multiplication of odd-degree elements is zero,
and the product 
$(x,\d w)\in E^{\mr{even}} X \x  E^{\mr{odd}} X \lt  E^{\mr{odd}} X$
descends from the multiplication of $\E W$
in the sense that
$
x\.\d w = \d\big(j_U^* i_U^*(x)\. w\big).
$
\eprop

\bpf
The additive isomorphisms follow from the reduction of the
Mayer--Vietoris sequence to
\[
0 \to E^{\mr{even}} X \os{i}\lt \E U \x \E V \lt \E W \os \d\lt E^{\mr{odd}} X \os{i}\to 0.
\]
The multiplication in the even subring follows
because $i$ is the ring homomorphism induced by $U \amalg V \lt X$.
The product of odd-degree elements $x,y \in E^{\mr{odd}}X$ is zero
by \Cref{thm:MV} since $\d$ is surjective.\footnote{\ 
	Alternatively, 
	since $i$ is injective on $E^{\mr{even}}X$ and vanishes on $E^{\mr{odd}}X$,
	we have $i(xy) = ix\.iy = 0$ so $xy = 0$.
}
To multiply an even-degree element $x$ with an odd-degree element $\d w$,
note that $\d$ is an $\E \mspace{-1mu}X$-module homomorphism
by \Cref{thm:MV},
so particularly $x\.\d w = \d (x\. w)$.
Now recall the module structure on $\E W$ is given by restriction
as $x\.w = (i_U \o j_U)^*(x)\.w$.
\epf

%
%
%

\brmk
In this paper, of course, we take $\E = \KG$.
In our previous joint work~\cite{CGHM2018},
we took $\E$ to be Borel cohomology $X \lmt H\Q^* (EG \ox_G X)$,
so that $\E(G/\G) = H\Q^*B\G$ is concentrated in even degree by Borel's theorem;
generally, given a nonequivariant cohomology theory $e^*$ 
such that $e^*(*)$ is torsion in odd degrees,
one could rationalize and take 
$\E$ to be rational Borel $G$-equivariant $e$-cohomology $e\Q_G^*$
so that $E^n(G/\G) = e\Q^n B\G$.
Since we have rationalized~\cite[Cor.~7.12]{rudyak},
the {\AHSS}s of CW-skeleta $B_n\G$ collapse at 
$E_2 = \H(B_n\G;\Q) \ox e^*(*)$, 
which is concentrated in even degree,
so that $\E(G/\G) = e\Q^*_G B\G$ is concentrated in even degree as well
and \Cref{thm:oddzero} applies.
The author is unsure how much demand there is for $e\Q_G^*$,
but has at least sighted the 
``Borel equivariant complex bordism''
functor $X\lmt MU_*(EG \ox_G X)$ in the wild.
\ermk